\documentclass[12pt]{article}
\usepackage{amsmath,amssymb,mathrsfs}
\usepackage{graphics}
\usepackage{graphicx}
\usepackage{titlesec}
\usepackage{subfigure}
\usepackage{boondox-cal}
\usepackage{color}
\usepackage{comment}
\graphicspath{{figures/png/}}
\usepackage{amsfonts}
\usepackage{latexsym}
\usepackage{url}
\usepackage{ifthen}  
\usepackage{bbm}
\usepackage{subfigure}

\newcommand{\ket}[1]{| #1 \>}
\renewcommand{\>}{\rangle}
\newcommand{\<}{\langle} 
\newcommand{\OUT}{|\mathrm{OUT}\>}
\newcommand{\IN}{|\mathrm{IN}\>}

\newcommand{\p}{\mathrm{prob}}
\newcommand{\JT}{\mathrm{JT}}

\newcommand{\thickbar}[1]{\mathbf{\bar{\text{$#1$}}}}
\newcommand{\ZZ}{\mathbf{Z}}
\newcommand{\Y}{\mathbf{Y}}
\newcommand{\X}{\mathbf{X}}
\titleformat{\paragraph}
{\normalfont\normalsize\bfseries}{\theparagraph}{1em}{}
\titlespacing*{\paragraph}
{0pt}{3.25ex plus 1ex minus .2ex}{1.5ex plus .2ex}
\begin{document}
\title{Jonckheere-Terpstra test for nonclassical error versus log-sensitivity relationship 
of quantum spin network controllers}
\author{E. Jonckheere,  S. Schirmer, and F. Langbein}
\maketitle

\abstract{Selective information transfer in spin ring networks 
by energy landscape shaping control has the property that 
the error $1-\p$, where $\p$ is the transfer success probability, 
and the sensitivity of the error to spin coupling uncertainties are %``positively correlated," 
%meaning that both are 
statistically increasing 
across a family of controllers of increasing error.   
The need for a statistical Hypothesis Testing of a concordant trend is made necessary by the 
noisy behavior of the sensitivity versus the error 
as a consequence of the optimization of the controllers in a challenging error landscape. 
%This already points to a defiance of the classical limitations.  
Here, we examine the concordant trend between the error and another measure of 
performance---the logarithmic sensitivity---used in robust control to formulate a well known fundamental limitation. 
%However, robust control rather advertises the logarithmic sensitivity  
%as the correct measure to formulate the fundamental limitations on achievable performance. 
Contrary to error versus sensitivity, the error versus logarithmic sensitivity 
trend is less obvious, because of the amplification of the noise due to the logarithmic normalization.    
%The problem is made complicated by 
%because its trend is corrupted by the noisy behavior of the logarithmic sensitivity 
%across controllers of increasing error numerically optimized in a challenging landscape.  
This results in the Kendall $\tau$ test for rank correlation between the error 
and the log sensitivity to be somewhat pessimistic with marginal significance level. 
Here it is shown that the Jonckheere-Terpstra test, 
because it tests the Alternative Hypothesis of an ordering of the medians of some groups of 
log sensitivity data, alleviates this statistical problem.  
This identifies cases of  
concordant trend between the error and the logarithmic sensitivity, 
a highly anti-classical features that goes against the well know sensitivity versus complementary sensitivity limitation.   
%with higher confidence than the Kendall $\tau$.} 

\section{Introduction}
\begin{comment}
Energy landscape shaping is used to control information transfer in
quantum networks, such as spin-1/2 rings and chains~\cite{Edmond_IEEE_AC}. 
Fabrication uncertainties of such networks at
nano-scale make it hard to ensure precise coupling strengths between the
spins as well as precise energy landscapes. This makes robustness of the
controls very important. Here we study the relation between the error
probability of the transfer of a single excitation and the logarithmic
sensitivity to such uncertainties. Contrary to classical robust control
results, the high performing (low probability error) quantum energy
landscape controllers also exhibit low logarithmic sensitivity.
\end{comment}

\subsection{Classical robust control bedrock}

One of the tenets of classical linear Single Degree of Freedom (SDoF) multivariable control~\cite{Safonov_Laub_Hartmann} is that the two fundamental 
figures of merit---tracking error and logarithmic sensitivity to model uncertainty---are in conflict. 
The former is quantified by the sensitivity matrix 
$S=(I+L)^{-1}$ and the latter by the complementary sensitivity $T=L(I+L)^{-1}$, 
where $L(s)$ is the input loop matrix.  
Specifically, 
\[ e_{\mathrm{track}}(s)=S(s)r(s),\]
 where $r(s)$ is an extraneous reference and $e_{\mathrm{track}}(s)$ is the tracking error, 
as shown in Fig.~\ref{f:classical_error}.   
$T(s)$ appears in the logarithmic sensitivity of $S(s)$ as $S^{-1}(dS)=(dL)L^{-1} T$. 
The conflict is obvious from $S+T=I$. 
%The conflict of SDoF configuration can be overcome by going to the 2-Degree of Freedom (2DoF) 
%configuration~\cite{2_deg_freedom_controller,Robust_2DOF_MIMO}, 
%or other architectures. 
The SDoF limitation can be overcome by 
a 2-Degree of Freedom (2DoF) 
configuration, as already pointed out by Horowitz~\cite[Chap. 6]{Horowitz} and recently made explicit 
in~\cite{2_deg_freedom_controller,Robust_2DOF_MIMO}. 

\begin{figure}[t]
\centerline{\scalebox{0.5}{\includegraphics{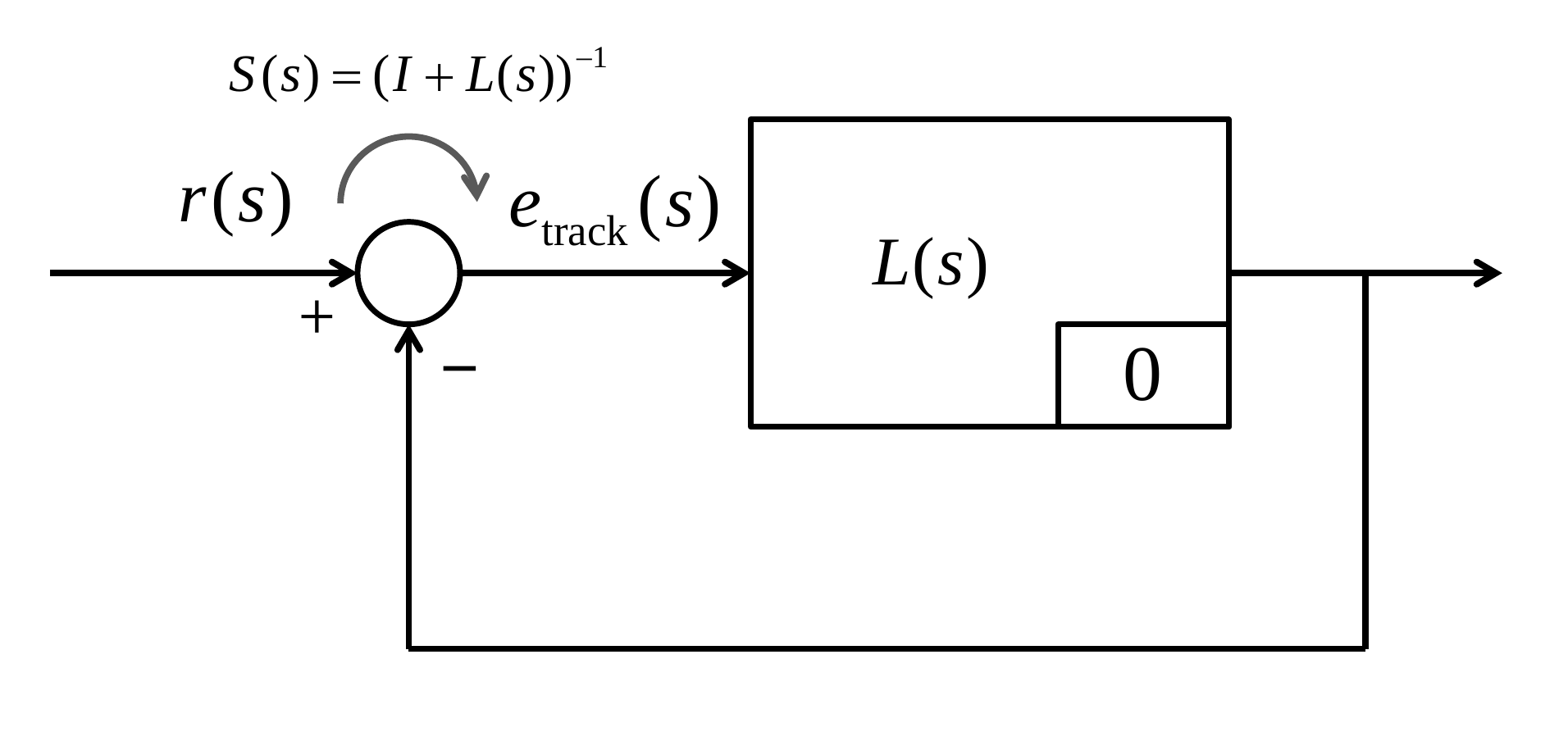}}}
\caption{Classical tracking error $e_{\mathrm{track}}$ in the classical control architecture}
\label{f:classical_error}
\end{figure}

\subsection{The quantum state overlap control problem}

This paper investigates whether this fundamental limitation survives in the quantum world, more precisely, in spin-$\tfrac{1}{2}$ networks 
where fabrication uncertainties at nano-scale make it hard to ensure precise coupling strengths between the spins. 
The answer is definitely negative 
for the class of quantum control problems where the objective is to achieve maximum overlap between the controlled wave function $|\Psi_D(t_f)\rangle$ at some final time $t_f$ and some reference or target wave function 
$|\Psi_{\mathrm{target}}\rangle$. In the preceding, $D$ denotes the controller. 
The controller belongs to a class of quantum control systems where the controller sole authority is to modify, 
{\it in a physically meaningful manner,} the parameters of the Hamiltonian~\cite{Tarn_SICON}.  
Even though the basic concepts developed here remain valid for the whole 
class of such controllers, here, however, we will more specifically consider 
controllers that induce energy level shifts~\cite{Chaos_soliton_fractals,quantum_rome,spin_network_curvature,chains_QINP,rings_QINP,time_optimal,Edmond_IEEE_AC}. 
More specifically in this paper,  
$D=\mbox{diag}\left( D_1, D_2,...,D_M\right)$ is a diagonal matrix of bias fields added to the Hamiltonian in some subspace.   
The concept is illustrated in Fig.~\ref{f:intro_picture}. 
In practice, we could envisage, e.g., electron spins in quantum dots whose energy levels can be controlled by voltages applied to surface gates~\cite{DiVincenzo}. 

\begin{figure}
\begin{center}
	\scalebox{0.4}{\includegraphics{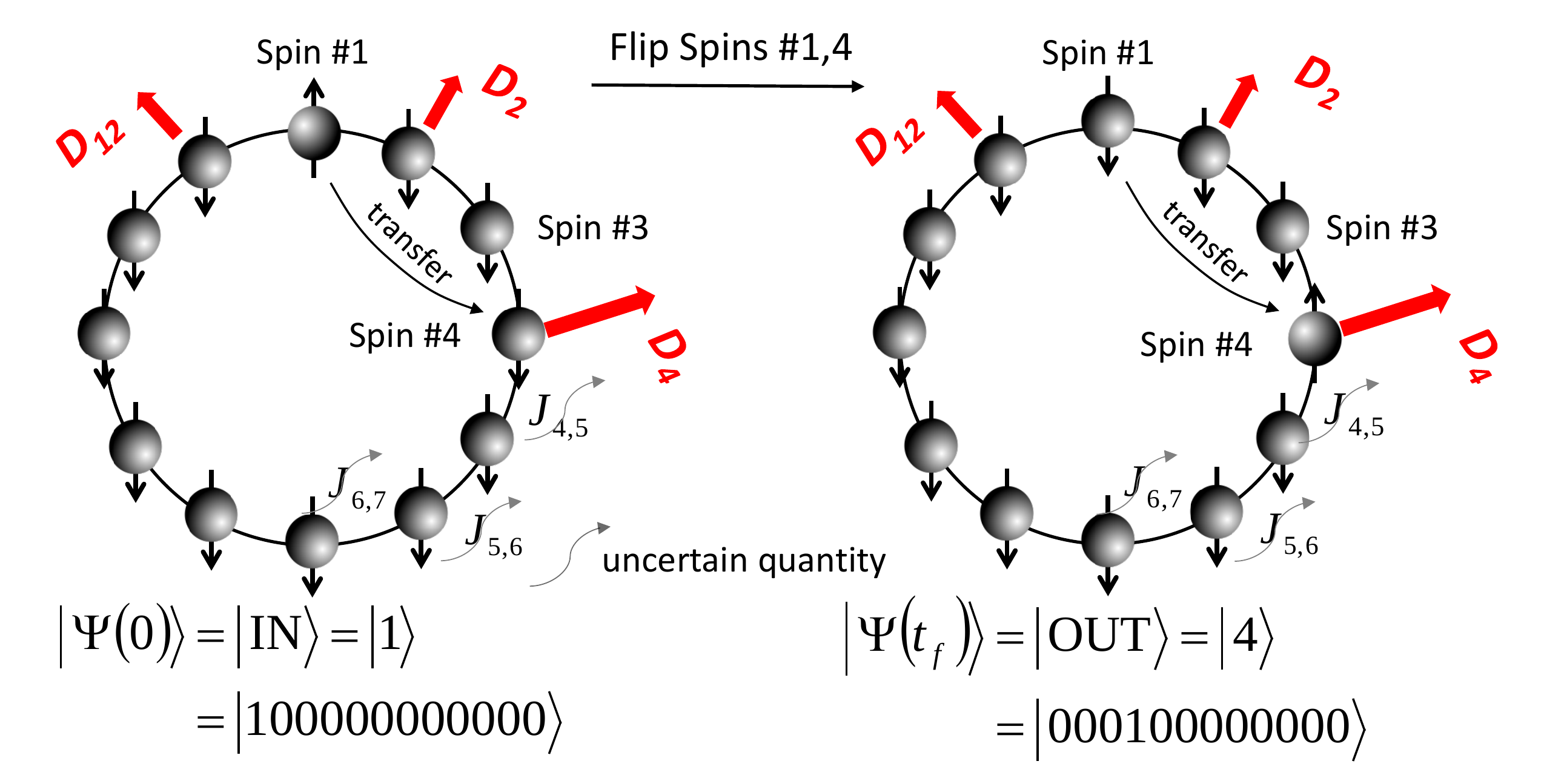}}
\end{center}
\caption{Illustration of the excitation transfer problem. 
Excitation (spin ``up") originally on 
spin \#1 has to be transferred to spin \#4. 
This is accomplished by bias fields $\{D_m\}_{m=1}^{12}$, 
despite uncertainties on the couplings $J_{m,m+1}$. 
(Not all bias fields are shown for clarity of the picture.) 
The initial (target) condition on the wave function $\Psi$ means that the 
INput (OUTput) excitation is on spin \#1 (spin \#4).}
\label{f:intro_picture}
\end{figure}

By ``maximum overlap," or ``maximum fidelity," we mean 
\begin{equation}\label{e:overlap} 
\max_{D,t_f} \left| \langle \Psi_{\mathrm{target}}|
\Psi_D(t_f) \rangle \right| \leq 1,    
\end{equation}
where $\langle \cdot | \cdot \rangle$ denotes the inner product in the complex Hilbert space. 
The upper bound on the fidelity is easily understood from the Cauchy-Schwartz inequality 
and the unit norm of the wave function. It turns out that this maximum overlap does not trickle down to the standard tracking problem
\[ \min_{D,t_f} \|\Psi_{\mathrm{target}}-\Psi_D(t_f)   \|, \]
but to a projective version of the tracking error~\cite{Edmond_IEEE_AC}:
\[ \min_{D,t_f,\phi} \|\Psi_{\mathrm{target}}-
e^{\jmath \phi}\Psi_D(t_f)   \|, \]
where $\phi$ is the global phase factor---a purely quantum mechanical concept. 
The classical-quantum discrepancy can now be understood as follows: Under some circumstances, the projective tracking error may be very small, yet the classical tracking error may be so large as to allow for small logarithmic sensitivity. 

In the same way as the conventional tracking error leads to a sensitivity matrix $S(s)$, 
the projective tracking error leads to an unconventional sensitivity matrix~\cite{Edmond_IEEE_AC}.  
The problem is that this sensitivity matrix does not easily lend itself to an analytical formulation of the limitations on achievable performance. 
Besides, the challenging landscape in which the optimization~\eqref{e:overlap} has to be conducted 
relative to both $D$ and $t_f$  
leads to ripples in the sensitivity versus error plot, calling for a statistical approach to determine concordance or discordance between the two figures of merit. 

\subsection{Statistical approach to classical-quantum discrepancy}

%The evidence that we propose in support of the classical-quantum discrepancy is statistical. 
We consider spin-rings as prototype 
quantum systems to test our hypothesis of anti-classicality. 
As shown in Fig.~\ref{f:intro_picture}, a $M$-spin ring is an assembly of $M$ spins arranged along a ring with near neighbor 
coupling $J_{m,m+1}$ of the XX or Heisenberg type. 
We chose spin rings because they are prototypes of quantum routers~\cite{rings_QINP} 
that have to transfer a single excitation (spin ``up") from one spin to another spin. 
To achieve this objective, time-invariant but spatially distributed bias fields $D_1$, $D_2$, 
..., $D_M$ are deployed. Such controller fits in the class of controllers mentioned earlier, as the $D_m$'s appear on the diagonal of the Hamiltonian. With this data, 
which includes the controller, the 
Schr\"odinger equation can be integrated analytically, 
and a closed form expression for the fidelity as well as its 
sensitivity can be derived~\cite{Edmond_IEEE_AC}. 

The problem is that there is no closed-form solution of the optimization~\eqref{e:overlap}. 
As shown in~\cite[Fig. 2]{time_optimal}, 
the error landscape is extremely challenging, and some optimization runs are successful at finding solutions very close to the upper bound, while the other solutions remain trapped in local minima with poor fidelity. 
The reward of having a variety of controllers achieving various levels of performance from optimal to poor is that it allows us 
to check---even quantify---the concordance or discordance between achievable transfer performance and  
the log sensitivity to the uncertain coupling parameters. 
This quantification is offered by the $Z$-statistic of the Kendall $\tau$ and more specifically 
the $|Z|$-statistic of Jonckheere-Terpstra on the data base~\cite{data_figshare} 
spanning across all rings from size 3 to 20, all transfers, and all coupling uncertainties. 

More specifically, the statistical analysis is done as follows:  The controllers \linebreak 
$\{D(n)\}_{n=1}^N$ from the database~\cite{data_figshare} are classified by increasing order of the error  \linebreak
$\sqrt{1-|\langle \Psi_{\mathrm{target}} | \Psi_{D(n)}(t_f)\rangle |}=:x_n$
they achieve. 
%Call this sequence $\{x_n\}$. 
Define $y_n$ to be the log-sensitivity of the $n$th controller. Classically, one would expect the sequences $\{x_n\}$ and $\{y_n\}$ to be discordant---that is, the error $x_n$ is  increasing while the sensitivity $y_n$ is decreasing with $n$. Contrary to classical control wisdom, in this quantum set-up the two sequences are concordant---that is, both $\{x_n\}$ and $\{y_n\}$ are increasing. 
Naturally, in this numerical set-up ``concordant" and ``discordant" have to be understood in a statistical sense. 
In general, the larger the $|Z|$-statistic of Jonckheere-Terpstra, 
the more the sequences are concordant.  

The purpose of this paper is to test the Null Hypothesis $H_0$ 
of no rank correlation between $\{x_n\}$ and $\{y_n\}$ using the Jonckheere-Terpstra $|Z|$-statistic  
on the variety of rings, 
subject the a variety of transfers, 
under the variety of parameter uncertainties compiled in the dataset~\cite{data_figshare}. 
In many cases, as identifed in Sec.~\ref{s:discussion}, $H_0$ is rejected in favor of the Alternative Hypothesis $H_A$ of concordance of  $\{x_n\}$ and $\{y_n\}$.

\subsection{Paper outline}

The paper is organized as follows: 
In Section~\ref{s:excitation_transport}, we review the spin network concept, the single excitation subspace, and we 
define the quantum excitation transport as the problem of having the solution to Schr\"odinger's equation 
move from an initial state of excitation to a target state of excitation. 
In Section~\ref{s:tracking_error_formulation}, 
the quantum excitation transport is contrasted with 
classical tracking control.
In Section~\ref{s:methodsI}, 
as an alternative to the analytical approach, 
we introduce the two statistical rank correlation tests---the Kendall $\tau$ and the Jonckheere-Terpstra tests---
which we propose to investigate whether 
the error and the logarithmic sensitivity are positively correlated. 
Section~\ref{s:methodsII} follows formal statistics methods and introduces the Type II error in the test. 
The statistical results specific to those controllers in the database~\cite{data_figshare} are presented in Section~\ref{s:resultsI},   
%, showing that in many cases 
%the Null Hypothesis of no rank correlation is rejected in favor of the Alternative Hypothesis 
%of a positive rank correlation, in contradiction with classical control. 
followed by Section~\ref{s:resultsII}, which  
shows that the power of the Jonckheere-Terpstra test as applied to the specific error versus sensitivity is within statistical gold standards. 
Finally, in Section~\ref{s:discussion}, we argue that, in excitation transport between nearby spins, classical limitations are overcome, 
while they tend to survive in case of transport between nearly diametrically opposed spins in rings. 
Appendix~\ref{a:related_tests} reviews some variants of the Jonckheere-Terpstra test, Appendix~\ref{a:left_tailed} 
reviews the left-tailed test, 
and Appendix~\ref{a:independence_formal}, the von Neumann (rank) ratio test, is presented as a test for independence of the observations.

\subsection{Notation}

The notation related to quantum physics is shown by the following table:
~\\

\begin{tabular}{r|l}
$\Psi$ & wave function of quantum ring\\
$H$ & Hamiltonian of uncontrolled quantum ring in single excitation subspace\\
$D$ & diagonal control Hamiltonian in single excitation subspace\\
$\IN$ & state of INput excitation into router\\
$\OUT$ & state of OUTput excitation out of router\\
$J_{m,m+1}$ & coupling between spins $m$ and $m+1$ \\
$M$ & number of spins in ring\\
$D_m(n)$ & $m$th component of bias field of $n$th optimization run\\
$N$ & number of fidelity optimization runs \\
$\mbox{err}$ & error $\sqrt{1-|\langle \Psi_{\mathrm{target}} | \Psi_D(t_f)\rangle |}$\\
$\mbox{prob}$ & probability of successful transfer $|\langle \Psi_{\mathrm{target}} | \Psi_D(t_f)\rangle |^2$\\
$\X_m,\Y_m,\ZZ_m$ & Pauli spin matrices of spin $m$ in network\\
$s^x,s^y,s^z$ & $2 \times 2$ Pauli matrices of single spin
\end{tabular}
~\\

The notation related to statistics is as follows:
~\\

\begin{tabular}{r|l}
$N$ & sample size of  (ring size, transfer, uncertainty) experiment\\
$x_n$ & independent variable (error 1-prob)) for sample $n$\\
$y_n$ & dependent variable ((log)sensitivity) for sample $n$\\
$I$ & number of bins in Jonckheere-Terpstra test\\
$i$ & a specific bin in Jonckheere-Terpstra test\\
$N_i$ & sample size in bin $i$\\
$\mu,\sigma^2$ & mean and variance, resp. \\
$s^2$ & unbiaised estimate of variance\\
$Z$ & normally distributed test statistic\\
$\JT$ & $=|Z|$, test statistic of Jonckheere-Terpstra\\
$p$ & $\int_u^\infty f_U(u)du$ where $f_U$ is the test statistic \\
$\alpha$ & significance level\\
$\mathrm{VN}$ & von Neumann ratio \\
$\mathrm{RVN}$ & rank von Neumann ratio  
\end{tabular}

\section{Excitation transport in networks of spins}
\label{s:excitation_transport}

\subsection{Network of spins}

The paper deals with the so-called single excitation, that is, a situation where one and only one spin in 
the network is ``up." It is however important to understand how this concept emerges 
from the general situation where as many as $M$ spins could be excited. 
The total Hamiltonian (including the controller) of a spin ring as the one shown in Fig.~\ref{f:intro_picture} is given by 
\begin{equation}\label{e:XYZmodel} 
\sum_{m=1}^{M} J_{m,m+1}(\X_m \X_{m+1}+\Y_m\Y_{m+1}+\varepsilon \ZZ_m\ZZ_{m+1})+\sum_{m=1}^M D_m\ZZ_m.
%+&J_{N_s,1}(X_{N_s} X_{1}+Y_{N_s}Y_{1} (+Z_{N_s}Z_{1})),
\end{equation}
In the above, $\X_m$, $\Y_m$, $\ZZ_m$ are the Pauli $x$, $y$, $z$ operators, respectively, of the spin $m$ in the ring,  
with the convention that 
$(\X,\Y,\ZZ)_{M+1}=(\X,\Y,\ZZ)_1$ to enforce the ring structure. More specifically,
\[ (\X,\Y,\ZZ)_m=I^{\otimes (m-1)}_{2 \times 2}\otimes s^{(x,y,z)} \otimes I^{\otimes (M-m)}_{2 \times 2},\]
where $s^{(x,y,z)}$ are the Pauli spin matrices
\[ s^x=\left(\begin{array}{cc}
0 & 1 \\
1 & 0
\end{array}\right), \; 
s^y=\left(\begin{array}{cc}
0 & -\jmath \\
\jmath & 0
\end{array}\right), \;
s^z=\left(\begin{array}{cc}
1 & 0 \\
0 & -1
\end{array}\right). \]
$J_{m,m+1}=J_{m+1,m}$ is the coupling strength between spins $m$ and $m+1$, with the convention that   
with $J_{M,M+1}=J_{M,1}$. Should $\varepsilon=0$ in the Hamiltonian, 
the ring is said to be XX, 
while it is said to be Heisenberg if $\varepsilon=1$. 
 
The operator $\ZZ=\tfrac{1}{2}\sum_{m=1}^{M} (I+\ZZ_m)$ counts the number of spins that are in the excited state. 
Since $\ZZ$ commutes with the Hamiltonian,   
the number of such spins remains invariant under the total motion. 
Define the single excitation subspace as the eigenspace of the +1 eigenvalue of $\ZZ$. 
%By definition, in the  
%single excitation subspace, the number of spins in the excited state remains exactly one. 
%In particular,  $\IN=e_i$ and $\OUT=e_o$,  
%where $\{e_m:m=1,...,M\}$ is the natural basis of $\mathbb{C}^{M}$ 
%and $1 \leq i,o \leq M$.  
%In the single excitation subspace, the dynamics reduces to~\eqref{e:single_excitation}. 
%``Single excitation" means that during the motion one spin and one spin only is ``up." 
%In particular,
In this single excitation subspace, the Hamiltonian of the $M$-ring  
reduces to the  $M \times M$ Hermitian matrix 
\begin{equation}
H+D=\left(\begin{array}{cccccc}
D_1 & J_{1,2} & 0 & \ldots & 0 & J_{1,M} \\
J_{1,2} & D_2 & J_{2,3} &  & 0 & 0 \\
0 & J_{2,3} & D_3 &        & 0 & 0\\
\vdots &&\ddots &\ddots&\ddots&                  \\
0 & 0 & 0 &        & D_{M-1} & J_{M-1,M} \\
J_{1,M} & 0 & 0 &  \ldots      & J_{M-1,M} & D_M
\end{array}\right)
\label{e:SES_Hamiltonian}
\end{equation}
and the wave function $\Psi \in \mathbb{C}^M$ is solution of the reduced Schr\"odinger equation  
\begin{equation}
\ket{\dot{\Psi}(t)}=-\jmath (H+D) \ket{\Psi(t)}, \quad \Psi(0)=\IN. 
\label{e:single_excitation}
\end{equation}
%
%to move from an initial single excitation state $\IN$ to some terminal single excitation state $\OUT$ 
%by means of some control $u(\cdot)$. 

\subsection{Quantum transport control}

In Eq.~\eqref{e:single_excitation}, 
$D$ could be time-invariant, time-varying, may or may not involve measurement feedback, 
the important point being that it should achieve high fidelity excitation transfer 
even in the presence of some uncertainties in $H$. 
By a well known, even fundamental, control paradigm that goes back to Bode, 
the latter can only be achieved if $D$ creates a feedback, possibly ``hidden," 
that wraps around the uncertainty. 
As already pointed out by Kosut~\cite{Kosut_Phy_Rev}, 
quantum gates can achieve both high fidelity operation and robustness with open-loop control, 
because the apparently ``open-loop" control creates a hidden feedback. 
This new paradigm can probably be best understood by splitting %(somewhat artificially) 
the open-loop $D$-controlled Schr\"odinger equation~\eqref{e:single_excitation}  
as a feedforward dynamics and a feedback control:
\begin{equation}
\begin{split}
\ket{\dot{\Psi}(t)}&=-\jmath H \ket{\Psi(t)} + u(t) , \quad \Psi(0)=\IN, \\
u(t) &= -\jmath D \ket{\Psi}(t).
\end{split}
\label{e:single_excitation_split}
\end{equation}
It clearly follows that, even when $D$ does not involve measurements as in~\cite{Edmond_IEEE_AC},  Eq.~\eqref{e:single_excitation} 
still involves some feedback that may be qualified as ``hidden" or better, ``field mediated," 
hence justifying the robustness properties. 
%Note, however, that in order to expose this hidden feedback 
%Eq.~\eqref{e:single_excitation} 

In the definition~\cite{Tarn_SICON} of quantum control as 
{\it ``tuning quantum interactions between matter and field, or field-field interaction"}, 
the $D$-controller would fall in the first category 
as, e.g., the electric fields from the gate electrodes control the energy
levels of electrons in quantum dots. 
%as the external bias magnetic field interacts with 
%the field developed by the spins. 
Here, as in~\cite{Tarn_SICON}, the $D$-controller is taken time-invariant, 
but spatially varying. 
In~\cite{Tarn_SICON}, the spatially varying controller is implemented by spatially modifying the dielectric constant of the medium of a wave guide. This has some commonality with 
our spatially varying bias field approach; 
however, our $D$-controller approach seems more related to the DiVincenzo architecture~\cite{DiVincenzo}. 

The drawback, however, of feed-backing $\ket{\Psi(t)}$ rather than the classical error 
$\Psi_{\mathrm{target}}-\ket{\Psi(t)}$ is that the controller has to be {\it selective}, that is, it must incorporate the 
knowledge of $\Psi_{\mathrm{target}}$, for otherwise the system has no way of knowing where to go. 

From the pure linear algebra viewpoint, 
observe that Eq.~\eqref{e:single_excitation_split} can be viewed as a linear feedback design, 
but a highly nonclassical one, as pointed out by Nijmeijer~\cite{bilinear_constant_input}. 
Indeed, the diagonal structure of $D$ 
takes the design outside the classical controllability pole placement problem and furthermore, because of the Hermitian property of $H$ and $D$, the poles can only be placed on the imaginary axis. 

\section{Tracking error formulation of quantum spin excitation transport}
\label{s:tracking_error_formulation}

 The excitation transport problem   
can, in some sense, be viewed as the problem of having $\ket{\Psi(t)}$ track $\OUT$. However, there are significant discrepancies between classical and quantum tracking control. 
First of all, the fundamental quantum figure of merit is not some error but 
the probability of successful transport of the excitation, 
or squared fidelity, $|\<\mathrm{OUT}|\Psi(t_f)\>|^2$, 
where $t_f$ is the time 
at which the excitation is read out. To simplify the exposition, assume that the probability achieves its maximum, 
$|\<\mathrm{OUT}|\Psi(t_f)\>|^2=1$, in which case it is easily seen that $\ket{\Psi(t_f)}=e^{-i\phi(t_f)}\OUT$, 
or equivalently 
$\ket{\Psi(t_f)}-e^{-i\phi(t_f)}\OUT=0$ for some global phase factor $\phi(t_f)$. 
More generally, it is not difficult to show that 
\begin{equation} 
\left\| \OUT - e^{\jmath\phi(t_f)} \ket{\Psi(t_f)} \right\|^2 = 2\underbrace{\left( 1-| \< \mathrm{OUT}|\Psi(t_f) \>|  \right)}_{\mathrm{err}^2(t_f)},    
\label{e:err}
\end{equation}
for
\[ \phi(t_f)=- \angle \< \mathrm{OUT}|\Psi(t_f)\>. \] 
It thus appears that the quantum transport problem of maximizing $| \< \mathrm{OUT}|\Psi(t_f) \>|$  
or its ``windowed" version $\frac{1}{\delta t}\int_{t_f-\delta t/2}^{t_f+\delta t/2}| \< \mathrm{OUT}|\Psi(t) \>|dt$ is equivalent to minimizing 
some ``tracking error''  
with the discrepancy that it is not required that the difference between the current state and the target state 
be small in the ordinary sense, but small in the sense of $\min_{\phi}\|\OUT -e^{i\phi}\Psi(t_f)\|$.  
The latter is related to the Fubini-Study metric~\cite[p. 31]{GriffithsHarris1994} 
on the complex projective space $\mathbb{C}\mathbb{P}^{M-1}$.  
%where $n_s$ is the number of spins (assuming single excitation).  

We will refer to the left-hand side of Eq.~\eqref{e:err} as the {\it projective tracking error.} 

\subsection{Classical-quantum controller structure discrepancies}

%Besides differences in the tracking error, the controller is not of the SDoF type, 
%but is in some sense a 2DoF controller. 
 
%is taken as a perturbation of the Hamiltonian~\eqref{e:SES_Hamiltonian}, 
The hidden feedback 
\begin{equation}
u(t)=-\jmath D\ket{\Psi(t)} 
\label{e:controller}
\end{equation}
formulation of the ``open-loop" controlled Schr\"odinger equation~\eqref{e:single_excitation}  
where $D$ is a diagonal matrix of spatially distributed biases 
makes the controller linear, 
as opposed to bilinear~\cite{Elliott,bilinear_constant_input}. 
Conceptually $D$ could still be time-varying, 
but here we focus on a time-invariant design. 
% that are shaping the energy landscape. 
%The time-invariance of $D$ 
As already said, in the latter, case the controller is linear time-invariant, but in the  
nonclassical sense of~\cite{bilinear_constant_input}. 
Yet another departure from classical control 
%in the sense that 
is that the controller is {\it selective,} that is, $D$ depends on both $\IN$ and $\OUT$. 
$\IN$ is the initial condition and, more importantly, $\OUT$ is to be interpreted as the reference. 
The controller is not driven by the tracking error, 
but depends on {\it both} the current state and the target state; 
from this point of view, the controller is of the 2DoF configuration. 
%linear in the state, as classically, 
%but the departure from the SDoF configuration is that 
%Note however that $D$ depends on $\OUT$ in some combinatorial fashion as there are only $M$ possible $\OUT$'s.  

Note that, because of the symmetry of the ring, $D(\IN,\OUT)$ depends only on the {\it distance} between 
$\IN$ and $\OUT$. 

Last but not least, 
the unitary evolution has the property that the controller is not asymptotically stable. 
Indeed, let $D$ be a controller that achieves $\| \OUT-e^{-\jmath(H+D)t}\IN \| \leq \epsilon$. 
Take an initial state $\IN'$ nearby $\IN$, that is, $\|\IN -\IN'\|=\eta$. Using the unitary property of the evolution and the triangle inequality, we derive 
\begin{eqnarray*}
\eta &=& \| e^{-\jmath(H+D)t}(\IN-\IN')\|   \\
&=&\|( e^{-\jmath(H+D)t} \IN - \OUT) +( \OUT-e^{-\jmath(H+D)t}\IN')\|   \\
&\leq& \epsilon + \|  \OUT-e^{-\jmath(H+D)t}\IN'    \| ,
\end{eqnarray*}
which yields
\[  \|  \OUT-e^{-\jmath(H+D)t}\IN'    \| \geq \eta-\epsilon .  \]
Thus, for an infinitesimally accurate controller ($\epsilon \downarrow 0$), 
the perturbed state will remain away from the target $\OUT$.  
% initial state $\IN'$,   
%that achieves $\IN \to \OUT$, for any $\IN' \ne \IN$, we will have $\|e^{-it(H+D)}\IN'-\OUT\|=1$, 
%$\forall t$.  
The latter has the consequence that the controller is not a classical asymptotically stabilizing controller;  
it is only Lyapunov stable. The latter property means that the state remains localized and won't diffuse. 
This means that in some cases~\cite{Chaos_soliton_fractals} our controller   
achieves 
{\it Anderson localization}~\cite{Anderson-58,short_to_Anderson,50_years}.  
%holding $D$ constant and symmetric has the consequence that the closed loop system  $\ket{\dot{\Psi}(t)}=-i (H+D) \ket{\Psi(t)}$ 
%is not asymptotically stable

With these significant departures from classicality, one wonders whether the fundamental 
error versus log sensitivity limitation 
is still in force. 
The problem is that the phase factor appearing in the quantum tracking error does not lead to  
a classical sensitivity function. 
In~\cite{Edmond_IEEE_AC}, a sensitivity matrix $\mathcal{S}(s)$  was defined 
via the Laplace transform $\hat{\mathcal{L}}$ of the projective tracking error as 
\begin{equation}
\hat{\mathcal{L}}(\OUT 1(t)-e^{i\phi(t)}\Psi(t)) = \mathcal{S}(s)\OUT 
\label{e:calS}
\end{equation}
for $\phi(t)$ achieving the minimum of 
$\|  \OUT-e^{i\phi(t)}\Psi(t)\|$. 
This sensitivity operator takes the form
\[ \mathcal{S}(s)= \frac{1}{s}I-\hat{\mathcal{L}}\left[e^{i \phi(t)}\right]\ast \left(sI+\jmath \left(H+D\right)\right)^{-1}P \]
where $P$ is a permutation such that $P\OUT=\IN$ and $*$ denotes  
the complex domain convolution~\cite{Edmond_IEEE_AC}. 
A clear relationship between $\mathcal{S}(s)$ and its sensitivity to parameters $J$ 
(a generic notation for $J_{k,k+1}$) in $H$ cannot be expected. 
For this reason, we propose a statistical approach based on a great many numerical optimization experiments. 

Note that, here, we define a sensitivity matrix without proceeding from a loop matrix as done classically. 
However, a fictitious loop matrix can be defined as 
\[ \mathcal{L}=\mathcal{S}^{-1}(I-\mathcal{S}) \]
and plugged in the feedback diagram of Fig.~\ref{f:projective_error}. Clearly, the conventional architecture is recovered, but for a very special loop matrix 
that embodies the projectivization of the error. 

\begin{figure}
\begin{center}
\scalebox{0.5}{\includegraphics{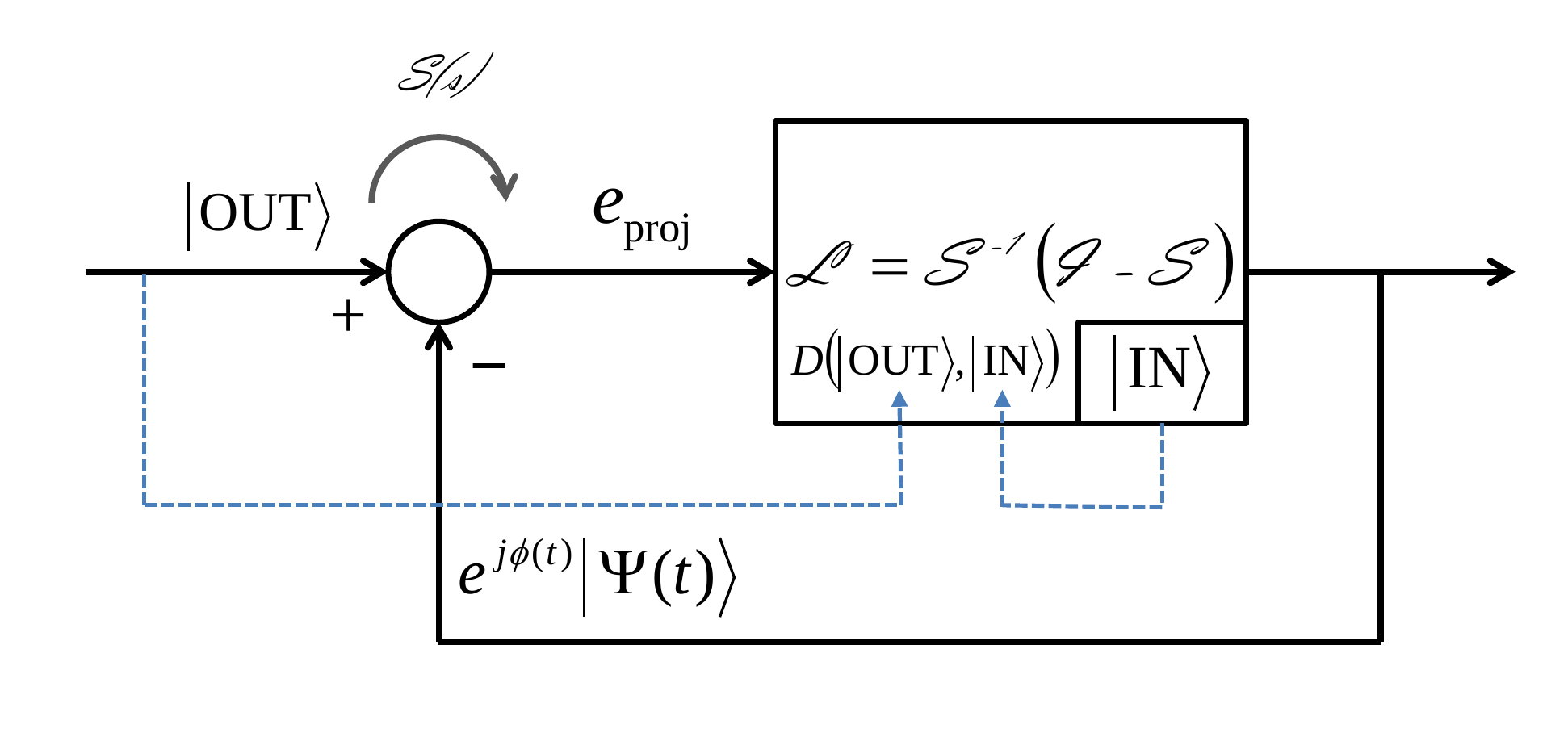}}
\end{center}
\caption{Projective error $\OUT-e^{j\phi(t)}\ket{\Psi(t)}$. 
Note that, contrary to the classical case of Fig.~\ref{f:classical_error}, the sensitivity matrix 
$\mathcal{S}$ is first defined, from which the fictitious loop function $\mathcal{L}$ is defined. 
The loop matrix is initialized with $\IN$. The dotted paths are nonclassical and indicate that the gain $D$ 
depends on both $\IN$ and $\OUT$. }
\label{f:projective_error}
\end{figure}

From~\eqref{e:calS}, it is clear that the sensitivity of the sensitivity $\mathcal{S}$ relative 
to $J$  amounts to sensitivity of $\mathrm{err}$ as defined by~\eqref{e:err}. From the classical control viewpoint, the log-sensitivity is 
\[ \frac{d \mathrm{err}}{dJ}\frac{1}{\mathrm{err}}=
-\frac{1}{4}\frac{1}{\sqrt{\p}}\frac{d \p}{dJ}\frac{1}{1-\sqrt{\p}}. \]
The right-hand side is easily derived from the definition of $\mathrm{err}$ taken from~\eqref{e:err}. 
If $\p \approx 1$, 
as the data base~\cite{data_figshare} retains only those controllers with an error not exceeding 0.1, 
then the above can be approximated as 
\[ \frac{d \mathrm{err}}{dJ}\frac{1}{\mathrm{err}}
\approx -\frac{1}{4} \frac{d \mathrm{prob}}{dJ}\frac{1}{(1-\sqrt{\mathrm{prob}})}\frac{2}{(1+\sqrt{\mathrm{prob}})}
= -\frac{1}{2} \frac{d\mathrm{prob}}{dJ} \frac{1}{1-\mathrm{prob}}\]
%then there is concordance between the left-hand side and the way the log-sensitivity is computed, 
%viz.,  
%$\frac{d\p}{dJ}\frac{1}{1-\p}=\frac{d\p}{dJ}\frac{1}{(1-\sqrt{\p})(1+\sqrt{\p})}
%\approx \frac{1}{2}\frac{d\p}{dJ}\frac{1}{(1-\sqrt{\p})}$.

\section{Methods---Type I error}
\label{s:methodsI}

%From a formal statistical viewpoint, this section deals with the Type I error, that is, the error of rejecting the  
%Null Hypothesis $H_0$ of no trend when it holds true. 

\subsection{Overview}

%\subsection{Utilization of results data base}

Here, as a first step towards an understanding of the error versus sensitivity issue, 
we proceed numerically 
by comparing 
the error $1-|\<\mathrm{OUT}|\Psi(t)\>|^2$  
and its (logarithmic) sensitivity to modeling uncertainties in $H$  
across a variety of controllers 
with error not exceeding 0.1 (see~\cite{data_figshare} for the data). % ordered by increasing error. 
Precisely, we considered {\it all} rings from $M=3$ to $M=20$ spins 
together with {\it all} transfers between any two spins. 
However, by symmetry, we can restrict ourselves to $\IN=\ket{1}$. 
The study therefore amounts to a total of 108 case-studies, where a case-study is defined by a number of spins $M$ and an $\OUT$ spin.  
For every $M\in[2,3,\ldots, 19,20]$ and every $(\IN=1,\OUT\leq \lceil M/2\rceil)$ pair, 
controllers $D$ were computed by numerical optimization runs of 
$1-|\<\mathrm{OUT}|e^{-i (H+D)t_f}|\mathrm{IN}\>|^2$ relative to $D,t_f$,  
either at the precise time $t_f$ or over a window around $t_f$, 
and controllers were ordered by increasing error, 
as explained in~\cite{time_optimal,Edmond_IEEE_AC} 
and as illustrated in Fig.~\ref{f:sensitivity_vs_logsensitivity}.  
Given a case-study $(M,\OUT)$ out of a total number of 108 case-studies, 
the number $N$ of time-windowed optimization runs, 
or controllers, were between 114 and 1998, with an average of 939 controllers. 

Note that we do not have a sampling of the set of {\it all} controllers. 
The set of controllers is the subset of those locally optimal controllers computed by the search algorithm 
and achieving an error not exceeding 0.1. 
%given by the search argorithm.}
%providing a statistically large sample.  

\begin{figure}
%\captionsetup{width=\textwidth}
\begin{center}
%\scalebox{0.5}{\includegraphics{sensitivity_vs_logsensitivity.png}}
\scalebox{0.5}{\includegraphics{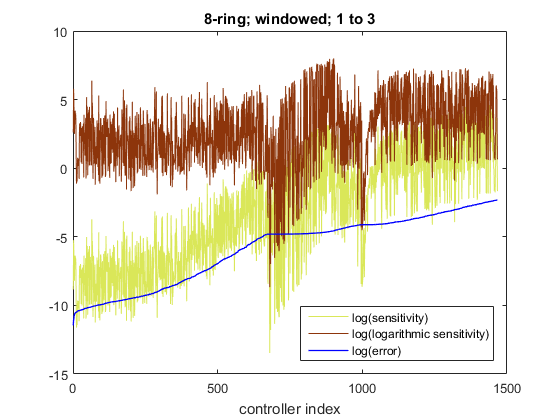}}
%\end{center}
\caption{Sensitivity versus logarithmic sensitivity. While the sensitivity is increasing with the error
with a Kendall $\tau$ of 0.6153, the behavior of the log-sensitivity is less trivial; 
nevertheless the Jonckheere-Terpstra test rejects the hypothesis of nonincrease of the log-sensitivity.}
\label{f:sensitivity_vs_logsensitivity} % label IMMEDIATELY after caption
\end{center}
%\label{f:sensitivity_vs_logsensitivity}
\end{figure}

The major difficulty is that 
the challenging error landscape and the potential for the solution to be trapped in some local minimums   
make the (absolute and logarithmic) sensitivity  versus error plots quite noisy,   
as shown by Fig.~\ref{f:sensitivity_vs_logsensitivity}, where controllers are ordered by increasing error.  
Despite this noisy behavior, the graph of Fig.~\ref{f:sensitivity_vs_logsensitivity} suggests a positive correlation between 
the sensitivity and the error (for this particular example). 
%It is pretty obvious from the plot that the error and the sensitivity are ``positively correlated." 
This observation is consistent with classical control; it is indeed   
easily seen that $dS=-S(dL)S$, meaning that if the error vanishes ($S=0$) so does the sensitivity ($dS=0$). 
However, 
the correlation between %the log sensitivity and the error
the logarithmic sensitivity $\left|\frac{d\p}{dJ}\frac{1}{1-\p}\right|$ 
relative to $J$-coupling uncertainties in $H$   
and the error  is not so obvious.  
%as easily seen from Fig.~\ref{f:sensitivity_vs_logsensitivity}.    
In order to make an {\it objective} statement about whether 
the logarithmic sensitivity versus error plot is increasing, decreasing, or inconclusive,    
we used two rank correlation test statistics: the Kendall $\tau$~\cite{Kendall}
and the Jonckheere-Terpstra statistic~\cite{Jonckheere_Terpstra,Terpstra_Jonckheere}.

%Since this latter relationship is less trivial and noisier than the absolute sensitivity case, 
%here we implement the nonparametric Jonckheere-Terpstra test, 
%which is believed to be more robust than the Kendall tau. 

%\subsection{More precise connection with classical tracking}

\subsection{Kendall $\tau$}

Given a set of independent, dependent variables pairs $\{(x_n,y_n)\}_{n=1}^N$, 
where $\{x_n\}_{n=1}^N$, $\{y_n\}_{n=1}^N$ are samples of random variables $\mathbf{x}$, $\mathbf{y}$, resp., 
the (estimate of the) Kendall $\tau$ is 
\[ \tau=\frac{\mbox{number of concordant pairs}-\mbox{number of discordant pairs}}{N(N-1)/2} \in [-1,1], \]
where a concordant pair is typically $((x_k<x_\ell) , (y_k<y_\ell))$ 
and a discordant pair is $((x_k<x_\ell),(y_k>y_\ell))$. 
The preceding assumes that there are no ties~\cite{Kendall_tau_with_ties}. 
A Kendall tau in $(0,1]$ means that the plot of $y$ versus $x$ is increasing---in the control context 
where $\mathbf{x}$ is the error and $\mathbf{y}$ the sensitivity, 
small (large) error implies small (large) sensitivity, a bit against traditional control wisdom. 

The mean and variance of Kendall $\tau$ are, respectively~\cite{tremendously_simplified},
\[ \mu_\tau=0, \quad \sigma_\tau^2=\frac{2(2N+5)}{9N(N-1)}.  \]
For large data set, the $\tau$ statistic
\[ Z_\tau=\frac{\tau}{\sigma_\tau} \]
is approximately normal, from which a test of significance can be drawn~\cite{Kendall_tau_significance}. 

A crucial condition is that the samples $\{y_n\}_{n=1}^N$ of $\mathbf{y}$ must be independent. 
This assumption can be justified by the randomness of the numerical optimizer running in an extremely 
complicated error landscape. In case of ``persistent data," there is a tendency towards an inflated value of the 
variance of $\tau$~\cite{Kendall_tau_for_persistent_data}. 

For the error versus sensitivity averaged over a small interval around $t_f$, 
the average Kendall $\tau$ over {\it all} rings from 3 to 20 spins and {\it all} transfers is 0.4535, 
indicating positive correlation, with a standard deviation of 0.2113, 
with an average $p$ of 0.001115741.  
However, for the logarithmic sensitivity, we obtained the less convincing values $\mu(\tau)=0.1925$  and $\sigma_\tau=0.2503$, with an average $p$ of 0.338925.  

The issue with the Kendall $\tau$ is that, 
when  it comes to the data $y$ increasing in an oscillatory fashion under an increase of $x$, 
the Kendall $\tau$ will find quite a few discordant pairs, even when on an average $y$ is obviously increasing. 
One remedy would be to smooth over $y$ and rerun the Kendall $\tau$ with the smoothed data. 
This of course would lead to a $\tau$ depending on the way the $y$ data has been smoothed over. 
Here we propose a different solution. 
The range of values of $x$ is decomposed in a certain number of groups, or ``bins,"  
and a Kendall $\tau$ like counting is made between groups, but not inside groups. 
This removes some of the discordant pairs and lead to a better figure of merit. 
This is the gist of the (nonparametric) Jonckheere-Terpstra test 
as it is applied to the present robust control problem. 

\subsection{Jonckheere-Terpstra test}
\label{s:JTstatistic}

Consider an independent variable, here the error $x=1-\p$ where $\p$ is the transfer success probability,  
and a dependent variable, here the logarithmic sensitivity of the probability relative to coupling errors $y=\left|\frac{d\p}{dJ}\frac{1}{1-\p}\right|$, where $J$ is the near-neighbor spin coupling strength. 
(Note that $1-\p$ is not the same error as $\mathrm{err}$, but this does not matter as the Jonckheere-Terpstra test is nonparametric, that is, it does not depend on the values but on the ranking of such values.) 
We want to show that $y(x)$ is statistically an increasing function (``positive correlation" between $x$ and $y$.) The range of values of $x=1-\p$ is decomposed in a certain number of groups 
such that the independent variable increases along the groups. To be formal, consider a partitioning   
of the values of the independent variable 
\[ \{x_n\}_{n=1}^N=X_1 \sqcup X_2 \sqcup ... \sqcup X_I \]
such that $\forall x_k \in X_i$, $\forall x_\ell \in X_j$ with $i<j$, 
we have $x_k \leq x_\ell$ with at least one strict inequality. 
With this grouping of the values of the independent variable, we construct a grouping of the corresponding 
values of the dependent variable:
\[ \{y_n\}_{n=1}^N = Y_1 \cup Y_2  \cup ... \cup Y_I, \quad Y_i :=y(X_i). \]
In each group of dependent variables, we compute the {\it median} of the population:
\[ \widetilde{Y}_1,\:\widetilde{Y}_2,\:...,\:\widetilde{Y}_I. \]
In the Jonckheere-Terpstra test~\cite{Jonckheere_Terpstra,Terpstra_Jonckheere}, the Null Hypothesis is
\[ H_0: \widetilde{Y}_1=\widetilde{Y}_2=...=\widetilde{Y}_I \]
and the Alternative Hypothesis is
\[ H_A: \widetilde{Y}_1 \leq \widetilde{Y}_2 \leq ...\leq \widetilde{Y}_I, \quad \mbox{with at least a strict inequality.} \]
The Jonckheere-Terpstra is a test for the Alternative Hypothesis. It is robust and avoids the noise in the log sensitivity because it argues on the medians. (The difficult part, though, is how to group the values of $1-\p$.) 

The statistic is derived from a counting of the number of cases favorable to the increasing property of 
$y$ relative to $x$ (the number of concordant pairs in the Kendall tau language). Precisely, 
we start with the Mann-Whitney $U$-statistic associated with the pair $(i,j)$ of groups:
\begin{equation}
\label{e:Uij}
 U_{ij}=\sum_{k=1}^{N_i}\sum_{\ell=1}^{N_j}\Phi(Y_j(\ell)-Y_i(k)), \quad i<j, 
\end{equation}
where 
\[ \Phi(z)=\left\{\begin{array}{cccc}
1 & \mbox{if} & z>0 &\\
1/2 & \mbox{if} & z=0&\mbox{(ties are counted as 1/2}) \\
0 & \mbox{if} & z<0&
\end{array} \right. \]
and $N_i=|Y_i|$ and $Y_i(k)$ denotes the $k$th element in $Y_i$. 
Defining $U=\sum_{i<j}U_{ij}$, the Jonckheere-Terpstra (JT) standardized test statistic 
\footnote{Note that in the original Jonckheere paper~\cite[Eq. 1]{Jonckheere_Terpstra} 
the statistic is rather defined, in our notation, as $U=\sum_{j=i+1}U_{ij}$. 
The $U=\sum_{i<j}U_{ij}$ is the formulation of the 
original paper by Terpstra~\cite[Eq. 3.1]{Terpstra_Jonckheere}. 
The Matlab {\tt JTtrend} function~\cite{JTtrend} follows the latter formulation.}
is $\mathrm{JT}=|Z|$, where 
\[ Z=\frac{U-\mu_U}{\sigma_U}. \]
Assuming that there are no ties~\cite{Jonckheere-Terpstra_with_ties}, 
the mean and the variance are, respectively~\cite{URS_in_JT_paper}, 
\begin{equation}\label{e:variance}
\begin{split}
\mu_U & =   \frac{N^2-\sum_{i=1}^I N_i^2}{4},  \\
\sigma_U^2 &= \frac{N^2(2N+3)-\sum_{i=1}^I N_i^2(2N_i+3)}{72},
\end{split}
\end{equation}
where $N=\sum_{i=1}^I N_i$. For a large data set, $Z$ is approximately normally distributed,  
from which the one-tailed $p$-value is computed as 
\begin{eqnarray} \label{e:pofZ}
p&=& \int_{u}^\infty f_U(u)du, \quad (u>0) \nonumber\\
&=&\frac{1}{2} \left( \int_z^\infty f_Z(z)dz +\int_{-\infty}^{-z} f_Z(z)dz \right), \quad (z>0) \nonumber\\
&=& \frac{1}{2} \left( 1-\mathrm{erf}\left(\frac{z}{\sqrt{2}}\right)  \right), \quad (z>0) \nonumber\\
&=&1-\frac{1}{2}\mathrm{erfc}\left(-\frac{\mathrm{JT}}{\sqrt{2}}\right), \quad (\JT=|z|). 
\end{eqnarray}
The various steps to compute $U$, $Z$, and $p$ from the $\{Y_i\}_{i=1}^I$ data  
are implemented in the Matlab function {\tt JTtrend.m} (see~\cite{JTtrend}), 
of which we have borrowed the notation.  

The Null Hypothesis $H_0$ of no trend is rejected if $p<\alpha$, where $\alpha$ is the significance level (by default 0.05) and accepted if $p>\alpha$. Equivalently, $H_0$ is rejected if 
$\JT > \JT_{\alpha}$ and accepted otherwise. The critical value $\JT_{\alpha=0.05}\approx 1.6557$ 
is easily verified from~\eqref{e:pofZ}.  If the Null Hypothesis holds with $f_Z$ normally distributed, $\alpha=0.05$ is the probability 
of wrongfully rejecting $H_0$; this is the Type I error. 

\noindent{\bf Remark:} There are other tests revolving around other ways to define the $U$-statistic. 
These are quickly reviewed in Appendix~\ref{a:related_tests}. 

\noindent{\bf Remark:}
Note that in the case the classical limitations are likely to hold, the left-tailed 
Jonckheere-Terpstra test should be implemented; see Appendix~\ref{a:left_tailed}.

\subsection{Domain of validity of Jonckheere-Terpstra test}
There are some conditions for the Jonckheere-Terpstra test to be applicable:
\begin{enumerate}
\item {\bf Independence of observations:} 
The requirement is that for each \linebreak 
$(M,\OUT)$ case-study 
the log-sensitivity ``observations" $\{y_n\}$ should be independent within each group and across all groups. 
This is empirically justified in 
Section~\ref{s:independence_intuitive}, where an argument based 
on the dynamics of the optimization algorithm that generates the data $\{y_n\}$ is developed.  
However, from the austere viewpoint of the mere data $\{y_n\}$ without reference as to how they were generated, 
this is the issue of securing some randomness in a series of observations, 
a problem that goes back to von Neumann~\cite{original_von_Neumann}. 
This is relegated to  Appendix~\ref{a:independence_formal}. 
\item {\bf Same group distribution shape:} The distributions of observations in each group must have the same shape and variability. This allows the Jonckheere-Terpstra test to be a test on the medians. 
Naturally, since we hope to find a trend in the data, this cannot hold true without some preprocessing, 
 typically the removal of the mean. For example, looking at Fig.~\ref{f:data-dt-11-1}, it is clear that removing the means  
over appropriate windows (``bins") will give equally distributed log-sensitivity data across the many windows (``bins"). 
The situation is somewhat more complicated in Fig.~\ref{f:data-dt-5-2}, where it is critical to correctly place the bins 
to secure consistent probability densities across the bins. 
Formally, the empirical cumulative distributions of the $\{y_{n_i}-\bar{y}_i\}$ data inside all bins  
should be compared. (See Section~\ref{s:empirical}.)
%In general, the larger the $I$, 
%the closer we get to the desired situation. 

\end{enumerate}

\subsection{Combining test statistics from independent experiments}

Let $\left(\JT_\eta,p_\eta\right)$ be the test statistic of the experiment $\eta$, that is, 
the $\JT$-statistic and the $p$ value across a variety of controllers 
for a given number $M$ of spins and a given transfer $\IN \to \OUT$. 
(For Jonckheere-Terpstra, we would include the number of bins $I$ in the experiment data.) 
$\JT_\eta$ and $p_\eta$ already allow for a Hypothesis Testing 
(accept or reject classical limitations) 
{\it for the given experimental set-up.} 
However, we want to do a Hypothesis Testing that transcends the particular experimental set up 
where the number of spins and the transfer are fixed. 
We want to do a Hypothesis Testing that spans across the many ring sizes and the many possible transfers. 
If the various experiments were $Z$-scored with $Z$ normally distributed (as in the Kendall $\tau$), 
the correct way to combine the various $Z_\eta$'s 
would be Stouffer's method~\cite{Stouffer_test_2013,Stouffer,Stouffer_test},   
$Z=(1/\sqrt{E})\sum_{\eta=1}^E Z_\eta$, 
since the resulting $Z$ is normally distributed, from which the $p$ is easily computed via~\eqref{e:pofZ}. 
A more intuitive way is Liptak's test~\cite{Stouffer_test_2013}, $Z=(1/E)\sum_{\eta=1}^E Z_\eta$. 
(The Fisher test~\cite{evol_biology}, $p=-2\sum_{\eta=1}^E \log p_\eta$,  directly combines the $p$'s.) 
The problem is that the Jonckheere-Terpstra test is $|Z|$-scored, not $Z$-scored. 
Here, we somewhat heuristically follow the Liptak method of just averaging the 
$|Z|$'s and the $p$'s and verifying from the data of Table~\ref{t:averages} that the relation~\eqref{e:pofZ},  
with $\JT$ and $p$ replaced by their means, holds up to several decimals. 
For the deployment of the Stouffer method, the reader is referred to~\cite{soneil_mu}. 

%Here we use the Liptak test and relegate the Stouffer test to~\cite{soneil_mu}. 

%Note that the Liptak test is more conservative than the Stouffer test. Indeed, to reject classical limitations, 
%we need a ``big” Z and even though the Liptak test already rejects the classical limitations the Stouffer test would reject them with a higher confidence since $Z_{\mathrm{Liptak}} < Z_{\mathrm{Stouffer}}$.

\section{Methods---Type II error and power of test}
\label{s:methodsII}

%Formally, a Type II error is the error of accepting the Null Hypothesis $H_0$ when it is untrue. 

The $U$ and the $|Z|$-statistics of Section~\ref{s:JTstatistic} hold true 
under the Null Hypothesis $H_0$ of no trend. 
Therefore, $\alpha$ as set to 0.05 is the probability of making the Type I error 
of rejecting $H_0$ when it holds true. 
\begin{comment}
Define the critical $|Z|_\alpha$ as 
%
\begin{equation}
\label{e:Zalpha}
 \alpha=1-\frac{1}{2}\mathrm{erfc}\left( -\frac{|Z|_\alpha}{\sqrt{2}} \right) . 
\end{equation}
%
\end{comment}
If we admit that this Type I error is small enough, 
%of rejecting $H_0$ when it is true, 
under the rejection of $H_0$, %with small Type I error if it is true, 
%the latter circumstances  
we would favor {\it some}  Alternative Hypothesis $H_A$. 
Instead of admitting the normal distribution 
$f_Z(z)=\frac{1}{\sqrt{2\pi}} \exp(-z^2/2)$ as true under $H_0$, 
we now shift the distribution to the right, $f_Z\left(z-\mu_{|Z|}^A\right)$, $\mu_{|Z|}^A >0$,  
and admit that the latter holds true under  
the Alternative Hypothesis. 
We then define the Type II error, $\beta_\alpha$, as the probability 
of failure to reject the Null Hypothesis when the specific Alternative Hypothesis holds:
\begin{equation}
\label{e:beta}
 \beta_\alpha\left(\mu_{|Z|}^A\right)=\int_{-\infty}^{Z_\alpha} f_Z\left(z-\mu_{|Z|}^A\right)dz
=1-\frac{1}{2}\mathrm{erfc}\left( \frac{Z_\alpha-\mu^A_{|Z|}}{\sqrt{2}}\right) . 
\end{equation} 

The power of the test given the Type I error $\alpha$ is defined as $1-\beta_\alpha\left(\mu_{|Z|}^A\right)$. 
Naturally, it depends on the admissible Type I error $\alpha$ and the conjectured $\mu^A_{|Z|}$ defining $H_A$, 
but it also depends on the overall sample size $N$ and the population $N_i$ of the bins. 
Since the $\mu^A_{|Z|}$ is guessed but otherwise a priori unknown,  
it is essential to assess the power of the test over all realistic $\mu^A_{|Z|}$'s. 
Usually, we require the power of a test to be 80\%.

\section{Results---Type I error} 
\label{s:resultsI}

\subsection{Conditions for test to be applicable}

\subsubsection{Independence of observations}
\label{s:independence_intuitive}

Here the observations are essentially the many (log)sensitivities achieved by the $D$-controllers 
obtained by running the optimization algorithm in the error landscape. 
%As already said, for the Jonckheere-Terpstra test to be applicable, observations need to be ``independent." 
We provide an intuitive justification of the independence of the observations, 
relegating the formal statistical argument to Appendix~\ref{a:independence_formal}. 

 The initial values chosen for the optimization of the fidelity relative to $D$ 
are a random sampling of
the domain of controllers. In most cases, the difference between the initial $D$-value and 
the maximum fidelity $D$-solution is small.  
%(would need to check that rigorously, but we have that
%data in the data set, I think).. 
Therefore, it appears that the random sampling should
result into a random sampling over the attraction domains for the
optimization algorithm. 
However, if the size of the domain of attraction is notably larger than the mesh of the random sampling of the space of controllers, 
then there is a bias. 
%Now, if these have notably different sizes, then there is a
%bias. 
If not, then the random initial value sampling matches a random
sampling of the maximum fidelity locally optimal controllers.  
Numerical experiments seem to indicate 
% ************Frank, check on this***************
that the same controller is not found twice over the runs. 
So this
would indicate that the attraction domains are much smaller than the
sampling density and that a random sampling of the space of initial controllers should lead to a random sampling of the resulting sensitivity/error data. 
  
It should be noted that contrary to time-series the ``time-stamp" $n$ assigned to be observations is not really a time. 
This can be explained by the way the results were derived.  
2000 independent optimization tasks were created for each  ``case-study," 
defined by a number of spins $M$ and a target spin $\ket{\mathrm{OUT}}$.
These were sent to a cluster and executed in parallel, such that the
sequence in which the results came out actually purely depends on the
cluster cores and the scheduler used. 
%This says nothing about the actual
%data, but about the sequence in which the tasks were executed (and they
Tasks were all mixed across all problems, with some rerun if a machine went down, etc. 
Each individual task selected an independent initial bias diagonal $D$-controller  
and initial time 
%(initial time values were taken from the time maxima of
%the corresponding chain evolution - that had to be done to make the
%optimization over time and bias feasible) 
according to a uniform distribution. 
So, the initial values were iid (using Matlab's
pseudo-random number generator) and as already argued this should lead to random (log)sensitivity results.  

\subsubsection{Empirical cumulative distribution}
\label{s:empirical}

Computation of the empirical cumulative distributions of the $\{y_{n}-\bar{y}_i\}_{n=1}^{N_i}$ data inside every bin after removal of the bin mean 
$\bar{y}_i=(1/N_i)\sum_{n=1}^{N_i} y_n$  reveal that they are fairly consistent. Sample results are shown in Fig.~\ref{f:ecdf}, 
where the partition of the data is uniform in $(N/I)$-observation groups, except for the last group. 
In geneal, the empirical cumulative distributions can be made closer by more careful grouping. 

\begin{figure}
\begin{center}
\subfigure{\scalebox{0.45}{\includegraphics{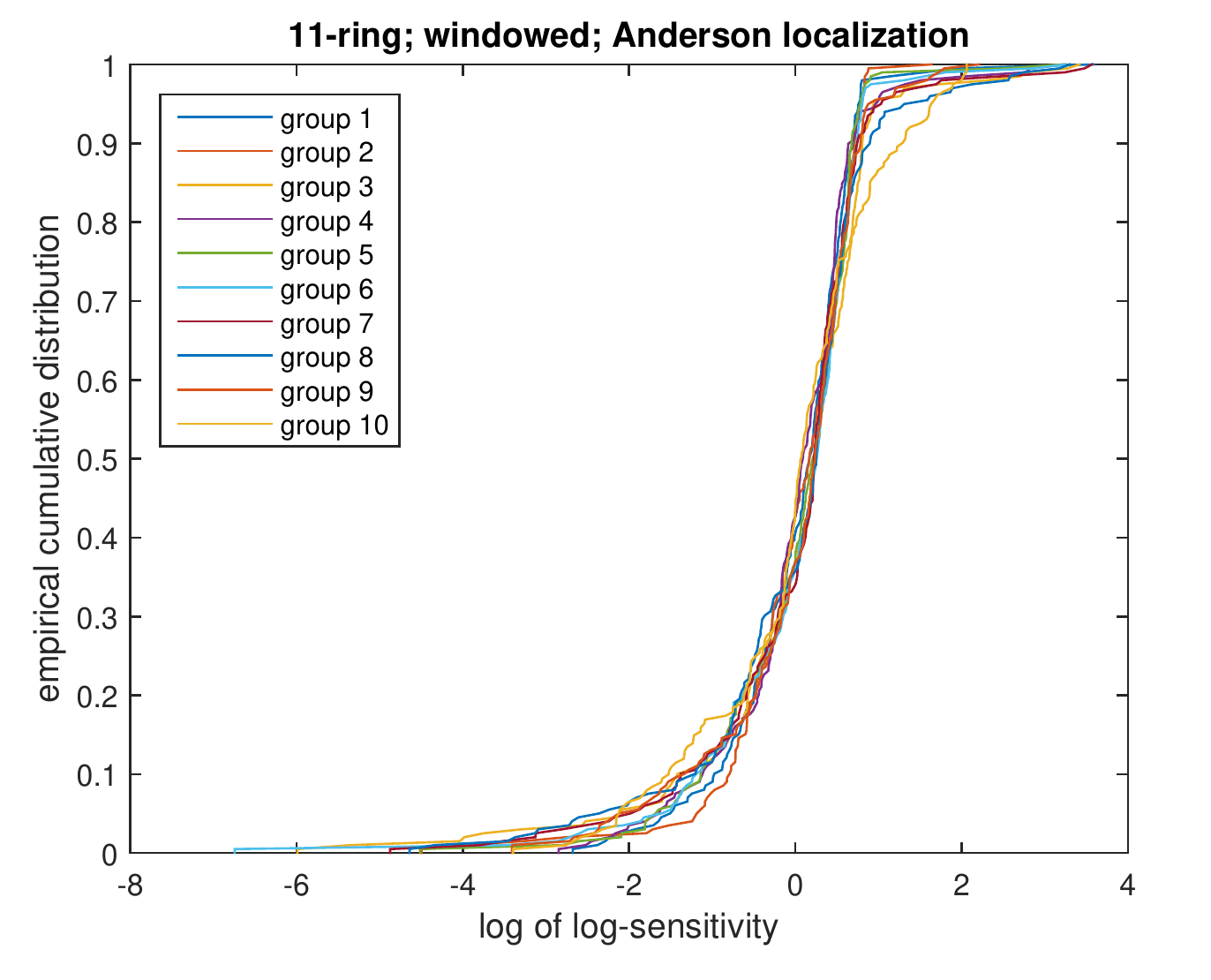}}} 
\subfigure{\scalebox{0.45}{\includegraphics{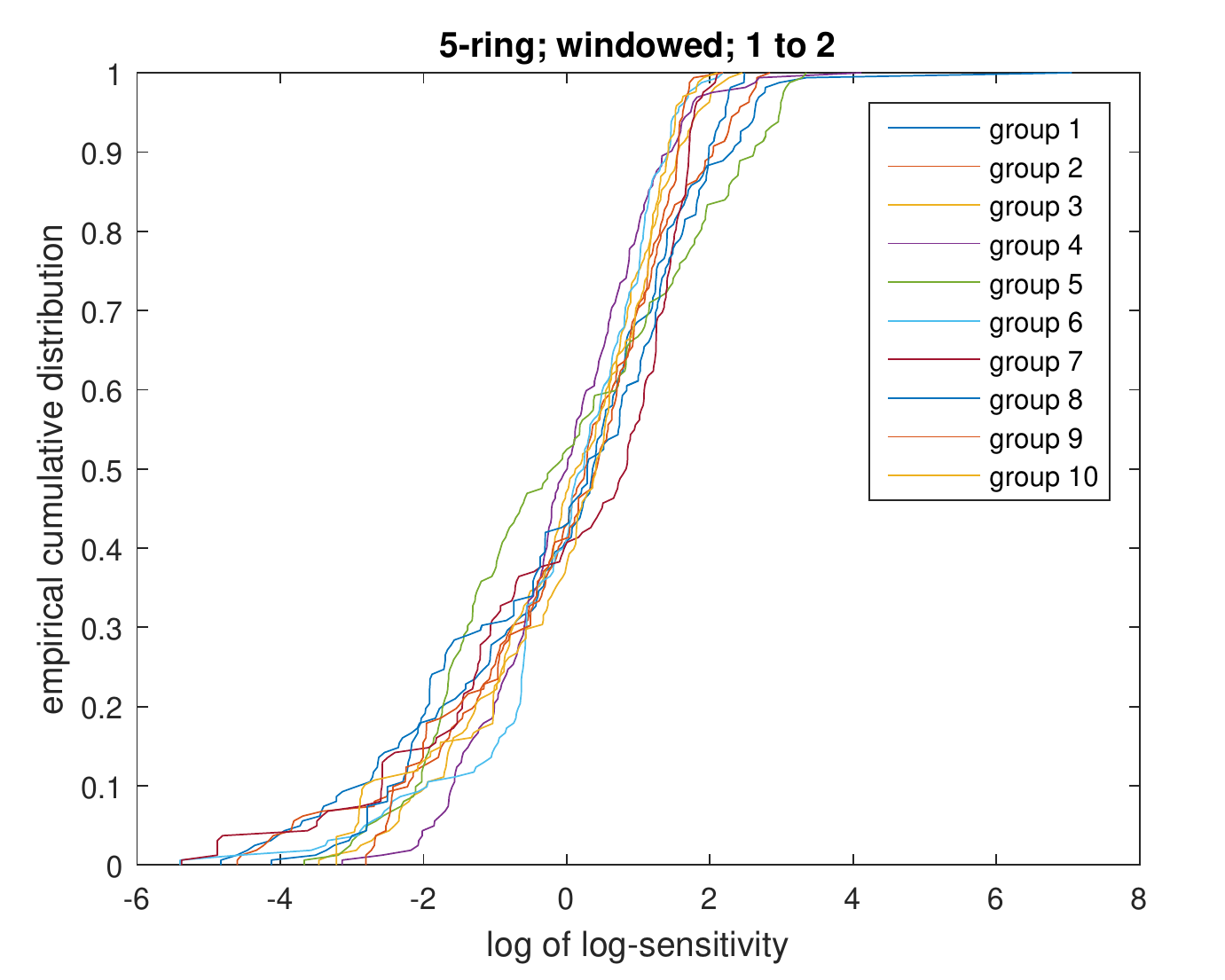}}}
\end{center}
\caption{Consistency of empirical cumulative distributions of observations divided in 10 groups of log-sensitivity data. 
Left: 11-ring with localization at $\ket{1}$; right: 5-ring with $\ket{1} \to \ket{2}$ transfer. 
The data before grouping is shown in Fig.~\ref{f:data-dt-11-1} and Fig.~\ref{f:data-dt-5-2}, resp.}
\label{f:ecdf}
\end{figure}

\subsection{Statistical analysis of error versus log-sensitivity relation}
\label{s:conflict_or_no_conflict}

Here we consider all case-studies of rings with $M=3$ to $M=20$ spins, 
with transport $\IN=\ket{1} \to \OUT$, with $\OUT$ ranging from $\ket{1}$ (Anderson localization) 
to $\lceil{\frac{M}{2}}\rceil$. 
This totals to an amount of 108 cases. 
By symmetry, this covers all cases of transfer of excitation 
from any spin to any other spin in networks of $M=3,4,\cdots,20$ spins. 

In each case-study among the 108 cases, we have $N$ pairs $\{x_n,y_n\}_{n=1}^{N}$.  
The $n$th value of the independent variable $x_n$ is the log of the error,  
$\log(1-\p_n)$,  
where $\p_n$ is the probability of successful $\IN \to \OUT$ transfer of the $n$th controller.   
(The log of the error allows for clearer graphing of the results yet it does not affect the ranking.) 
The errors are in increasing order $x_k \leq x_\ell$ for $k<\ell$. 
The dependent variable takes values
$$ y_n=y(x_n)=\frac{1}{2}\log\left(\sum_m\left|\frac{d\p_n}{dJ_{m,m+1}}\frac{1}{1-\p_n}\right|^2\right), $$
where $J_{m,m+1}$ is the $m\mbox{-} (m+1)$ spin coupling strength and the sum is extended over all couplings. 
In our data base, $N$ ranges 
%from 300 up to 1,500. 
from 114 up to 1998. %double-check this!!!!
The set of pairs is divided into $I$ groups, $\{(X_i,Y_i)\}_{i=1}^I$,   
where we took 
$I=3,10, 100$. 

For each data set $\{x_n,y_n\}_{n=1}^{N}$ corresponding to a certain number of spins 
and a certain $\IN \to \OUT$
transfer, the JT statistic $Z$ and the $p$ value were computed using 
the Matlab {\tt JTtrend} function~\cite{JTtrend}.  
From the $p$ value, a ``reject/accept" decision on the Null Hypothesis $H_0$ of no trend was taken 
consistently with a significance level $\alpha=0.05$. 
The average results over all case studies are shown in Table~\ref{t:averages}.
\begin{table}[t]
\caption{Analysis of the $\JT:=|Z|$ Jonckheere-Terpstra statistic over whole data base~\cite{data_figshare}   
(The $\min p=0$ is up to 4 decimals.)}
\begin{center}
\begin{tabular}{|c||ccc|c||ccc|c|}\hline
\rule{0pt}{2.5ex}
$I$ & $\min \JT$ & $\thickbar{\JT}$ & $\max \JT$ & $s_{\JT}$ & $\min p$ & $\bar{p}$ & $\max p$ & $s_p$ \\\hline\hline
3   &   0.0472 &  11.3055 &  33.0028 & 11.1729    &  0       & 0.0617  & 0.4812    & 0.1182     \\\hline
10  &   0.0433 &  12.1213  & 34.7227  & 12.1481    &  0       & 0.0595  & 0.4827    &  0.1156     \\\hline
100 &   0.364 &  12.2425 & 35.0389  & 12.3419    &  0       & 0.0630  & 0.4855    & 0.1235      \\\hline
\end{tabular}
\end{center}
\label{t:averages}
\end{table}

From Table~\ref{t:averages} the following conclusions can already be drawn:
\begin{enumerate}
\item There is not much difference between the $I=3$, 10 and 100 (number of bins) cases, 
except for the outlier $\min \mathrm{JT}=0.364$, $I=100$. The problem is that the data set contains 
controllers for a $(M,\OUT)$ case-study and that this case-study has a sample set of only $N=150$ controllers. 
Clearly, the arrangement of the $\{(x_n,y_n)\}_{n=1}^{150}$ data in ``bins" of 100 defeats the purpose of robustification of the results by arguing on the medians of the bins. 
For this particular case, it turns out that the two medians do not conform to the rest of the results. 
Note that there is another $N=114$ case-study, but for that one the 2 medians conform with the other results.   
\item The mean $p$-value is borderline between ``accept" $H_0$ (no disagreement with classical limitations) 
and ``reject" $H_0$ (disagreement with classical limitations), with a slight tipping of the balance toward ``reject." 
(Recall that $\alpha=0.05$.)  
\item $\min p = 0$ (up to 4 decimals)  means that there are cases in strong disagreement with classical limitations---the log sensitivity increases with the error.
\item $\max p \approx 0.48 \gg 0.05$ means that there are cases where there is not enough evidence to disagree with the classical limitations---meaning that the log sensitivity does not have trend relative to an increase error. 
\item Comparing Kendall $\tau$ with Jonckheere-Terpstra it is absolutely obvious that 
$\bar{p}_{\mathrm{Jonckheere-Terpstra}} \ll \bar{p}_{\mathrm{Kendall}~\tau}$. Clearly the Jonckheere-Terpstra 
test implicitly filters the oscillatory logarithmic sensitivity data and renders a result 
with significantly higher confidence than the Kendall $\tau$.  
\end{enumerate}

The upshot is that a simple relation like the classical $S+T=I$ cannot, in general, be expected in the quantum transport setup---except for the Anderson localization case, that is, holding a state of excitation at a single spin, or securing a successful ``transfer" $\ket{1} \to \ket{1}$. In this case indeed $p$ is consistently vanishing up to 4 decimals, rejecting the no trend hypothesis in the log sensitivity and pointing towards an increase of the log sensitivity with the error.   This anti-classical behavior is not surprising, as the Anderson localization is probably the quantum transport case that most significantly departs from classical concepts. 

\subsection{Case studies}

\subsubsection{Case-study: Anderson localization: ``reject'' classical limitation}

We consider the case of an 11-ring with the $\ket{1} \to \ket{1}$ ``transfer." Fig.~\ref{f:data-dt-11-1} shows that the various figures of merit are not conflicting---quite to the contrary, they are consistent. The detail of the experiment is shown in Table~\ref{t:data-dt-11-1}. Clearly, the ``reject" decision is consistent with the visual appearance 
of the log sensitivity plot. 

\begin{figure}[t]
\begin{center}
\scalebox{0.5}{\includegraphics{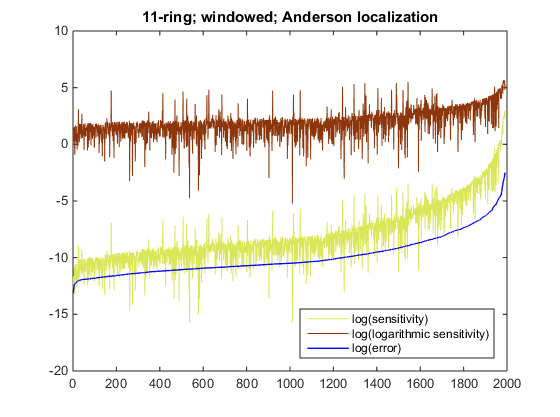}}
\caption{
Sensitivity, logarithmic sensitivity, and error plotted on a logarithmic scale (for better comparison) 
versus index $n$ of controller (2000 of them) for Anderson localization around spin 1 in an 11-ring.  
%Various nonconflicting figures of merit of an 11-ring under Anderson localization around spin 1. 
All 3 figures of merit are concordant, indicating that Anderson localization is anti-classical.
}
\label{f:data-dt-11-1} % label ALWAYS after caption and IMMEDIATELY after caption
\end{center}
\end{figure}

\begin{table}[h]
\caption{Details of the 11-ring Anderson localization experiment on Jonckheere-Terpstra test 
of Null Hypothesis of no trend between error and logarithmic sensitivity}
\begin{center}
\begin{tabular}{|c||c|ccc|}\hline
$I$ & Kendall tau & $|Z|$ & $p$ & Null Hypothesis \\\hline\hline
3      & 0.4483 & 26.5509 & 0 & ``rejected"    \\\hline
10    & 0.4483 & 29.5768 & 0 &  ``rejected"  \\\hline
100  &0.4483 & 29.8896 &  0 & ``rejected"     \\\hline
\end{tabular}
\end{center}
\label{t:data-dt-11-1}
\end{table}

\subsubsection{Case study: ``reject" classical limitation}

Anderson localization is not the only case where an anti-classical behavior is observed, as shown by 
the strongly increasing trend of the log sensitivity in the case of a 5-ring under $\ket{1} \to \ket{2}$ transport shown in 
Figure~\ref{f:data-dt-5-2}. The details of the analysis is shown in Table~\ref{t:data-dt-5-2}. 

\begin{figure}[t]
\begin{center}
\scalebox{0.5}{\includegraphics{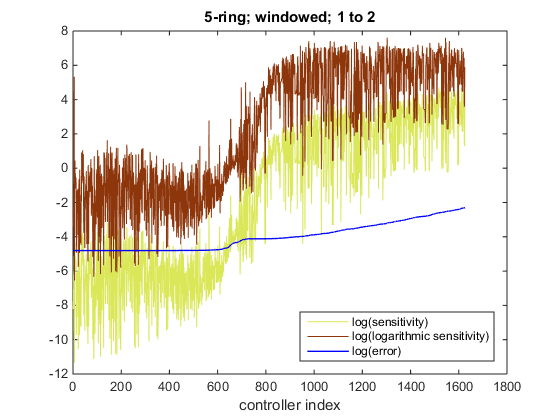}}
\caption{
Sensitivity, log-sensitivity, and error plotted on a logarithmic scale 
versus index $n$ of controller (1600 of them) for $\ket{1} \to \ket{2}$ transport in a 5-ring. 
Observe the strong concordance between the sensitivity and the log-sensitivity. 
Even though the concordance with the error is weaker, 
it still indicates anti-classical behavior.
%Strong increasing trend of the log sensitivity of a 5-ring under $\ket{1} \to \ket{2}$ transport
}
\label{f:data-dt-5-2} % label ALWAYS after caption and IMMEDIATELY after caption
\end{center}
\end{figure}

\begin{table}[t]
\caption{Details of the 5-ring under $\ket{1} \to \ket{2}$ transport experiment 
on Jonckheere-Terpstra test 
of Null Hypothesis of no trend between error and logarithmic sensitivity}
\begin{center}
\begin{tabular}{|c||c|ccc|}\hline
$I$ & Kendall tau & $|Z|$ & $p$ & Null Hypothesis \\\hline\hline
3      & 0.58 & 33.0028 & 0 & ``rejected"    \\\hline
10    & 0.58 & 34.7227 & 0 &  ``rejected"  \\\hline
100  &0.58 & 35.0389 &  0 & ``rejected"     \\\hline
\end{tabular}
\end{center}
\label{t:data-dt-5-2}
\end{table}

\subsubsection{Case-study: borderline ``accept/reject" classical limitation}

As ``borderline"  case, we choose a 14-ring with $\ket{1} \to \ket{6}$ transfer. The log sensitivity plot of Fig.~\ref{f:data-dt-14-6} shows first an 
increasing trend and then a decreasing trend relative to the error, 
which explains the mixed ``accept/reject" decision shown in Table~\ref{t:data-dt-14-6}. 

\begin{figure}[t]
\begin{center}
\scalebox{0.5}{\includegraphics{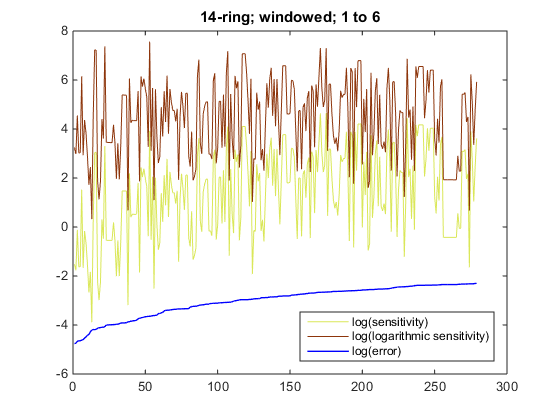}}
\caption{Sensitivity, log-sensitivity and error plotted on a logarithmic scale versus 
index $n$ of controller (260 of them)  
for a 14-ring under $\ket{1} \to \ket{6}$ transport. 
While the error and sensitivity are concordant, 
the error and log-sensitivity are marginally concordant, 
because of the decreasing trend of the log-sensitivity as of controller 200.
}
\label{f:data-dt-14-6} % label ALWAYS after caption and IMMEDIATELY after caption
\end{center}
\end{figure}

\begin{table}[t]
\caption{Details of the 14-ring $\ket{1} \to \ket{6}$ transport experiment 
on Jonckheere-Terpstra test 
of Null Hypothesis of no trend between error and logarithmic sensitivity 
(Recall that $\JT_{\alpha=0.05}\approx 1.6557$.)}
\begin{center}
\begin{tabular}{|c||c|ccc|}\hline
$I$ & Kendall tau & $|Z|$ & $p$ & Null Hypothesis \\\hline\hline
3      & 0.0575 & 1.4875  & 0.0684  & ``accept"    \\\hline
10    &  0.0575 & 1.5144  & 0.065 &  ``accept"  \\\hline
100  &  0.0575 & 1.6696  & 0.0475  & ``reject"     \\\hline
\end{tabular}
\end{center}
\label{t:data-dt-14-6}
\end{table}

\subsubsection{Case study: ``accept" classical limitation}

Here we consider one of the best illustrative case of no increase of the log sensitivity. 
We consider the case of a 15-ring with the $\ket{1} \to \ket{6}$ transfer. Fig.~\ref{f:data-dt-15-6} shows that the 
logarithmic sensitivity has no trend compared with the error, as confirmed by the details of Table~\ref{t:data-dt-11-1} and the ``admit" 
the Null Hypothesis decision. 

\begin{figure}[t]
\begin{center}
\scalebox{0.5}{\includegraphics{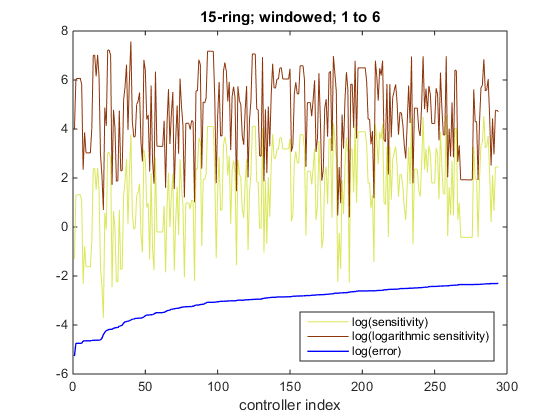}}
\caption{Sensitivity, log-sensitivity, and error plotted on a logarithmic scale 
versus index of controller (290 of them) in a 15-ring under $\ket{1} \to \ket{6}$ transport.
While the sensitivity and the error are concordant, 
the log-sensitivity shows no trend---even possibly a slightly decreasing trend---versus $n$ 
indicating rather classical behavior.
%The ``no trend" behavior of the logarithmic sensitivity relative to the error in a 15-ring under $\ket{1} \to \ket{6}$ transport.
}
\label{f:data-dt-15-6} % label ALWAYS after caption and IMMEDIATELY after caption
\end{center}
\end{figure}

\begin{table}[ht]
\caption{Details of the 15-ring $\ket{1} \to \ket{6}$ transport experiment 
on Jonckheere-Terpstra test 
of Null Hypothesis of no trend between error and logarithmic sensitivity}
\begin{center}
\begin{tabular}{|c||c|ccc|}\hline
$I$ & Kendall tau & $|Z|$ & $p$ & Null Hypothesis \\\hline\hline
3      & -0.0285 & 0.1789  &0.4290  & ``accept"    \\\hline
10    &  -0.0285 & 0.7549  & 0.2252 &  ``accept"  \\\hline
100  &  -0.0285 & 0.2052  & 0.4200  & ``accept"     \\\hline
\end{tabular}
\end{center}
\label{t:data-dt-15-6}
\end{table}

\section{Results---Type II error and power of test}
\label{s:resultsII}

We compute the power of the Jonckheere-Terpstra test for 
$\alpha$ in a neighborhood of the 0.05 significance level decided upon in the previous sections, 
and in the generic experimental situation where $M=1000$ and $I=10$. 
We set $\alpha$, compute the variance from~\eqref{e:variance}, 
compute the critical $|Z|_{\alpha}$ by setting $p=\alpha$ in Eq.~\eqref{e:pofZ}, 
and finally compute $\beta_\alpha$ from~\eqref{e:beta}.  
%Note that for $\alpha=0.05$, $|Z|_{\alpha=0.05}\approx 1.6446$.  
The results are shown in Fig.~\ref{f:multi_alpha_power}. 
In order to acheive the ``gold standard" of 80\% power of the test, we need to have $\mu_{|Z|}^A \geq 2.5$, 
which is certainly achieved for the cases shown in Tables~\ref{t:data-dt-11-1} and~\ref{t:data-dt-5-2} 
where trends are detected. 

We also considered the power of the Jonckheere-Terpstra test for various number $M$ of spins, for various number $I$ of ``bins," 
to conclude that the power is not visibly affected by those quantities.
\begin{figure}[t]
\begin{center}
%\scalebox{0.5}{\includegraphics{nice_power_cropped}}
%\scalebox{0.5}{\includegraphics{multi_alpha_power_cropped.pdf}}
\scalebox{0.5}{\includegraphics{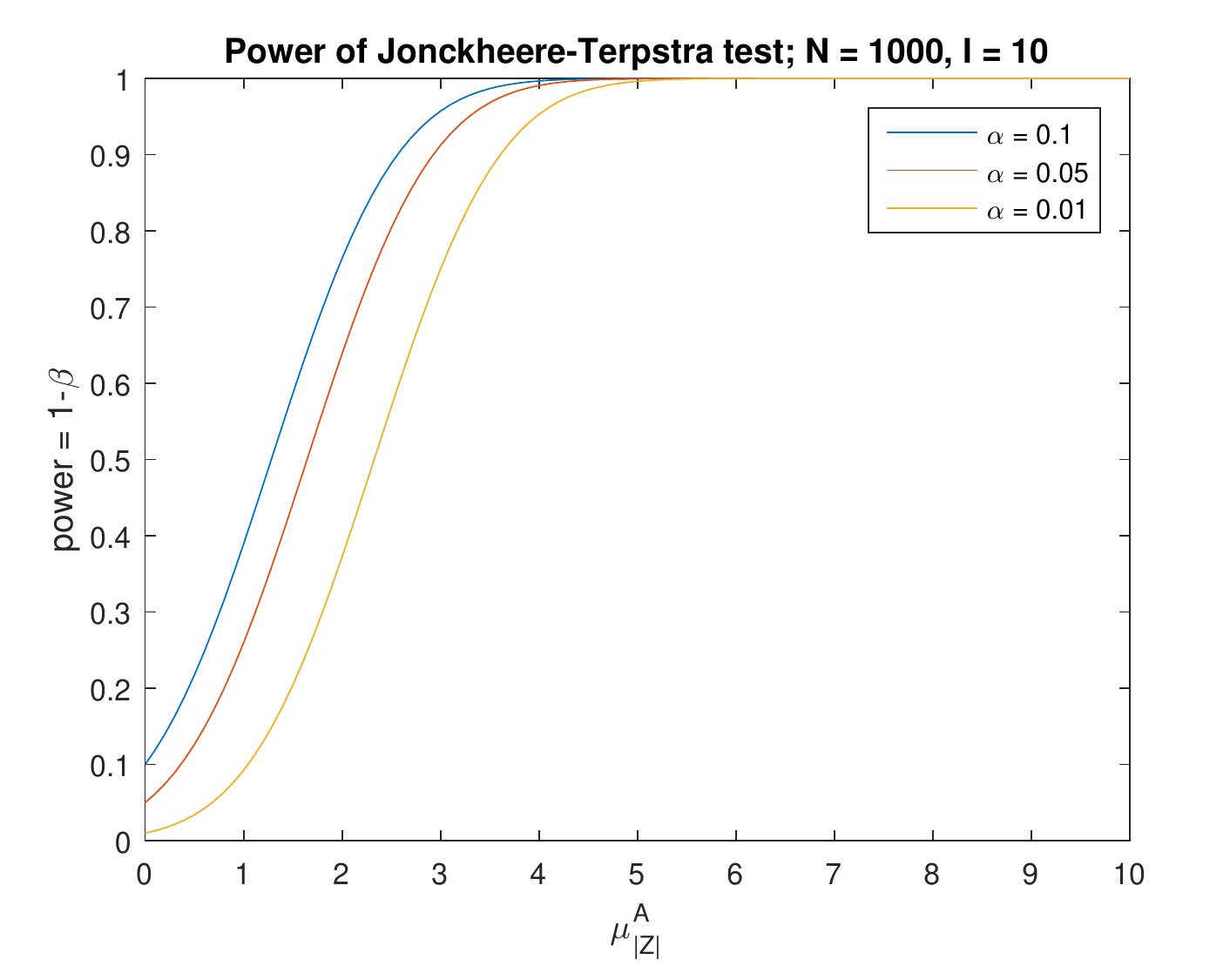}}
\end{center}
\caption{Power of Jonckheere-Terpstra test for various $\alpha$'s versus mean $\mu^A_{|Z|}$ under Alternative Hypothesis}
\label{f:multi_alpha_power}
\end{figure}

 %On the other hand, the power of the test depends on $\alpha$, as shown in Fig.~\ref{f:multi_alpha_power}. 

\section{Discussion: Dependency of error versus log-sensitivity relation on $(\IN,\OUT)$}
\label{s:discussion}

In the previous study, the data incorporated {\it all} cases, up to symmetry, of $\IN \to \OUT$ transfers, 
for all $M$ ranging from 3 to 20, 
with an overall positive concordant trend between error and log sensitivity. 
Here we examine how much the classical/anti-classical behavior depends on the 
relative position of the $\IN$ and $\OUT$  spins. 
The overall Jonckheere-Terpstra $|Z|$-data with $\IN=\ket{1}$   
of Section~\ref{s:conflict_or_no_conflict} is divided into three $\OUT$-groups ($I=3$) that roughly correspond to a decomposition of the right-half of the ring into 3 equally-sized sectors:
\begin{itemize}
\item $Y_1$: $\mathrm{JT}$-data of ($120^\circ<\mbox{angle}(\ket{1},\OUT)<180^\circ+\epsilon),$
\item $Y_2$: $\mathrm{JT}$-data of~~($60^\circ<\mbox{angle}(\ket{1},\OUT)<120^\circ$), 
\item $Y_3$: $\mathrm{JT}$-data of~~~~($0^\circ<\mbox{angle}(\ket{1},\OUT)<60^\circ$),  
\end{itemize}
as illustrated in Fig.~\ref{f:3_groups}. 
If $M$ is not divisible by 3, we arrange the $M$-spin data such that $|Y_3| \geq |Y_2| \geq |Y_1|$ with at least a strict inequality. The more specific grouping of the data is shown in the contingency  Table~\ref{t:contingency}.

\begin{figure}
\begin{center}
%\scalebox{0.6}{\includegraphics{3_groups_cropped.pdf}}
\scalebox{0.6}{\includegraphics{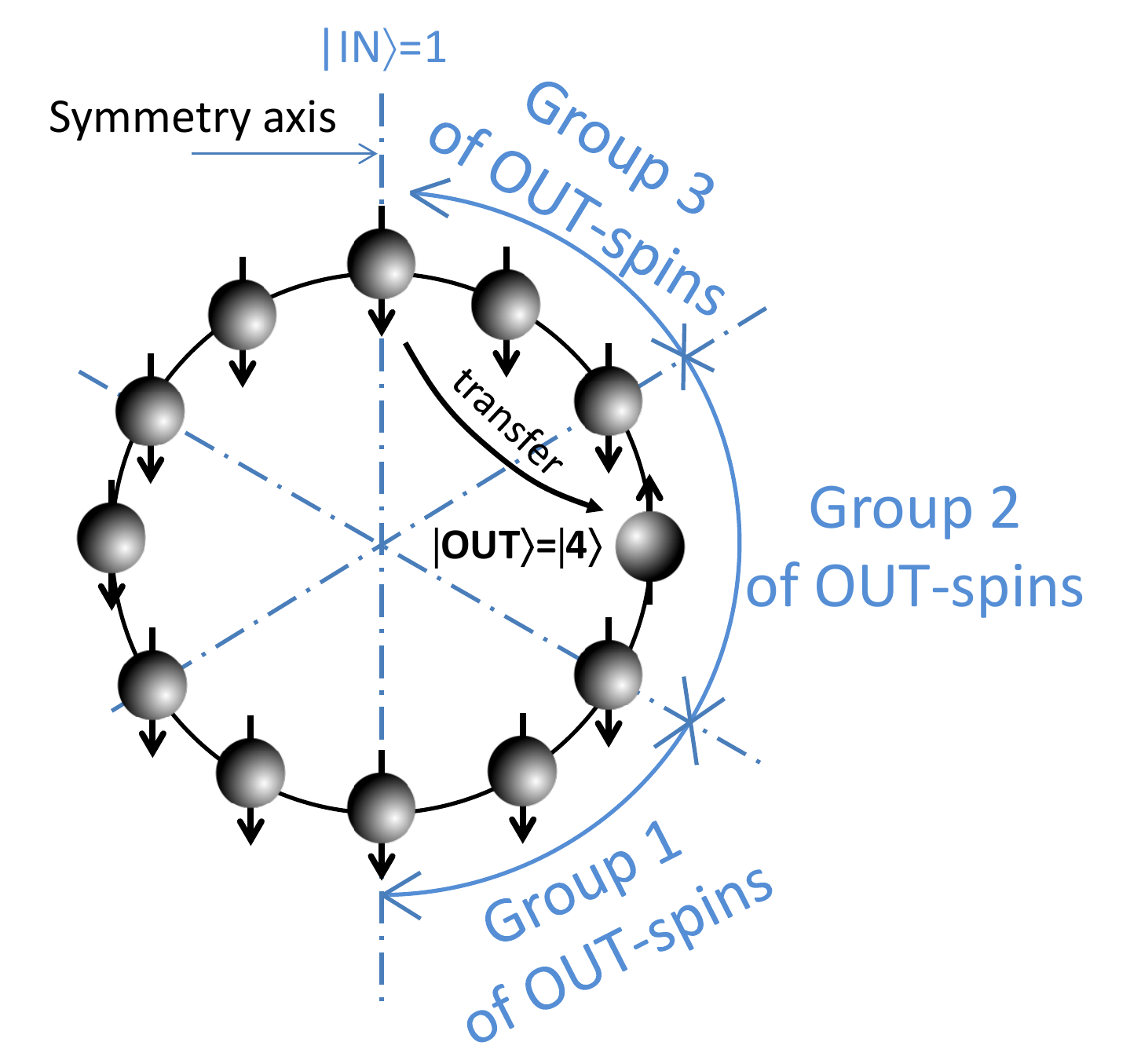}}
\end{center}
\caption{The 3 groups of $\OUT$-data to assess dependency classical/anti-classical behavior 
on $\OUT$-position around the ring}
\label{f:3_groups}
\end{figure}

\begin{table}[t]
\caption{Contingency table of data: $[M:M']-[\OUT:\OUT']$ denotes all $Z$-data, $I=3$, 
pertaining to a number of spin between $M$ and $M'$ with initial spin state $\IN=\ket{1}$ and target spin state ranging 
from  $\OUT$ to $\OUT'$.}
\begin{center} 
\begin{tabular}{|c|c|c|}\hline\hline
$Y_1$          & $Y_2$          & $Y_3$ \\\hline\hline
10-5           &  $[10:12]-[3:4]$ & $[10:12]-[1:2]$ \\
$[11:12]-[5:6]$  & $[13:14]-[4:5]$  & $[13:18]-[1:3]$ \\
$[13:14]-[6:7]$  & $[15:18]-[4:6]$  & $[19:20]-[1:4]$ \\
$[15:16]-[7:8]$  & $[19:20]-[5:7]$  & \\
$[17:18]-[7:9]$  &                & \\
$[19:20]-[8:10]$ &                &\\\hline
\end{tabular}
\end{center}
\label{t:contingency}
\end{table}

With the JT-data arranged in the $Y_1$, $Y_2$, and $Y_3$ bins, 
we examine whether there is a trend in the JT data across the three bins. 
The Jonckheere-Terpstra test rejects the Null Hypothesis of no trend with $|Z|=7.6283$ and $p=0.0000$ 
(up to 4 decimals) for the Alternative Hypothesis of a trend 
$\tilde{Y}_1 \leq \tilde{Y}_2 \leq \tilde{Y}_3$, with at least one strict inequality.  
Therefore, when the spins of excitation transfer $\IN \to \OUT$ are not too far apart, 
the design behaves anti-classically 
(error and log sensitivity increase together). When they become nearly anti-podal, 
then the design behaves classically with the conflict between error and log sensitivity. 

Note that this conclusion is supported {\it by the available dataset,} which   
contains only those controllers with a largest error of 0.1. 
%That's why the number is not consistent. 
This in particular means we have fewer controllers for the
longer distance transitions on the ring, as it was considerably harder to find these.   
Possibly the conclusion could be invalidated by better  optimizers 
able to find better controllers at long distance transport.

\section{Conclusion \& Future research directions}

As already observed in~\cite{Edmond_IEEE_AC}, the quantum transport problem can be re-formulated in the classical 
control setup only at the expense of a complicated sensitivity matrix $\mathcal{S}(s)$ 
that embodies the projectivity of the quantum tracking error.  
Given the projective sensitivity matrix $\mathcal{S}(s)$, the fundamental limitations, if any,  
are not easy to come by. Here we have developed a statistical approach based  
on a sample set~\cite{data_figshare} of numerically optimized controllers. 
%Having constructed 
Precisely, for {\it every} ring from 3 to 20 spins and {\it every} transfer on such ring,   
a fairly large data set of locally optimal controllers, 
arranged by increasing order of their transfer errors, was constructed~\cite{time_optimal,data_figshare}. 
With the controllers at hand, 
we investigated whether the logarithmic sensitivity anti-classically increases 
with the error using the Kendall tau and the Jonckheere-Terpstra tests, 
with a preference for the latter as it gives higher confidence. 
Out of all case-studies constructed from the whole dataset~\cite{data_figshare}, it appears that  
there are cases that clearly show anti-classical behavior ($H_0$ rejected), while others show
classical behavior ($H_0$ accepted). 
%the Jonckheere-Terpstra Null Hypothesis of no trend may or may not be rejected. 
The former---rejection of no trend in favor of an increasing trend---is a challenge to the classical limitations 
that say that the error and the logarithmic sensitivity should be in conflict. 
By a further analysis, it was shown that for transfers between nearby spins, 
the classical limitation does not hold, 
while it tends to be recovered for transfers between distant spins. 

The results derived here are based on a data set that retains, 
among other numerically optimized controllers, 
only those achieving 
a probability error no greater than 0.1. 
%Those controllers were confronted with small perturbations of the coupling constants  together with 
% a differential formulation of  fundamental limitations. 

The issue of large as opposed to differential parameter variations 
is addressed in~\cite{ssv_mu,soneil_mu}, 
where a structured singular value argument proves that challenge to the classical limitation remains in force. 
Note that~\cite{soneil_mu} not only considers coupling errors but also field focusing errors 
and that the same $\mu$-analysis argument~\cite{ssv_mu,soneil_mu} is able to cope with initial state preparation errors.   
But a more challenging robustness problem consists in evaluating the classical limitation 
in the context of the gap between a model like~\eqref{e:XYZmodel} and some real-life quantum components, 
like some Copper compounds~\cite{Heisenberg_Copper}, 
which approach the Heisenberg model~\eqref{e:XYZmodel}, but will never quite match the model. 
The same challenge applies to the DiVincenzo architecture~\cite{DiVincenzo} versus its models.

Finally, observe that all that precedes applies to coherent quantum dynamics. 
However, when the ring is subject to collective dephasing, the classical limitations 
tend to re-appear~\cite{classicality}.

%Allowing more controllers 
%in the data set along with {\it large} variations will be considered in a further paper, 
%but a preview that uses structured singular values is already available 
%in~\cite{ssv_mu}. 

\section*{Acknowledgment}

The authors wish to thank Professor Stanislav Minsker, Department of Mathematics, University of Southern California, for drawing our attention to the Type II error. Many thanks to Capt. Sean O'Neil, 
US Military Academy,  
for drawing our attention to the Stouffer test.  
Sophie Schirmer and Frank Langbein are
supported by the Ser Cymru NRN AEM grant 82.

\appendix

\section{Some related tests}
\label{a:related_tests}

There are many extensions/refinements of the Jonckheere-Terpstra test~\cite{powerful_nonparametric}. 
In case of small data sets, a modified version of~\eqref{e:Uij} is proposed as 
\[ U_{ij}=(j-i)\sum_{k=1}^{N_i}\sum_{\ell=1}^{N_j}\Phi(Y_j(\ell)-Y_i(k)), \quad i<j. \]
Another recently proposed version~\cite{powerful_nonparametric} is 
\[ U_{ij}=(r_{j\ell}-r_{ik})\sum_{k=1}^{N_i}\sum_{\ell=1}^{N_j}\Phi(Y_j(\ell)-Y_i(k)), \quad i<j, \]
where $r_{j\ell},r_{ik}$ denote the position (rank) of $Y_j(\ell), Y_i(k)$ in the combined data. Finally, 
yet another extension proposes a confidence interval~\cite{critique_of_JT}. 

It was observed that the first refinement of the Mann-Whitney $U$-statistic does not change the overall results and conclusion. 

\section{Left-tailed Jonckheere-Terpstra test}
\label{a:left_tailed}

In case the classical limitations are likely to hold, 
the Jonckheere-Terpstra test should be  organized around the Alternative Hypothesis 
%This test might be too conservative for our purposes, because the classical limitations of control theory 
%claim that equality
%%
%\[ \widetilde{Y}_1=\widetilde{Y}_2=...=\widetilde{Y}_I \]
%%
%cannot be achieved; only inequality 
%
\[ \widetilde{Y}_1 \geq \widetilde{Y}_2 \geq ...\geq \widetilde{Y}_I, \]
%
%can be achieved. So we need to tweak the Jonckheere-Terpstra test as a test to reject the above classical robust control ordering. 
that is, the log sensitivity is decreasing with increasing error. 
The test is analogous to the classical one, but in the opposite tail. 
Rejection of the Null Hypothesis in favor of the above Alternative Hypothesis is more likely to happen with the 
instantaneous performance optimizing controllers. 
This is left to a further paper.

\section{(Rank) von Neumann test}
\label{a:independence_formal}
A qualitative argument in favor of the randomness of the results of the search algorithm was presented in Sec.~\ref{s:independence_intuitive}, 
but a quantitative analysis stills needs to be set up. 

\subsection{von Neumann ratio test}

The genesis of the von Neumann test~\cite{original_von_Neumann,earlier_von_Neumann_ratio_test} is to decide whether a trend exists in a series of observation $\{y_n\}_{n=1}^N$ {\it totally} ordered by 
the variable $n$, 
usually thought to be the time. The von Neumann test of independence relies on the paradigm  that a trend compromises the randomness of a time-series. To quantify this observation, define the mean square successive difference
\[ \delta^2=\frac{1}{N-1}\sum_{n=1}^{N-1} (y_{n+1}-y_n)^2 \]
and the (slightly biased) variance estimate
\[ s^2=\frac{1}{N}\sum_{n=1}^N (y_{n}-\bar{y})^2, \quad \left(\bar{y}=\frac{1}{N}\sum_{n=1}^N y_n\right). \]
The von Neumann ratio is defined as 
\[ \mathrm{VN}=\frac{\delta^2}{s^2}. \]
Taking $y_n$ linear in $n$, hence giving $y_n$ a trend, it is easy to see that $\lim_{N \to \infty} \delta^2/s^2 \to 0$. 
Intuitively, small $\mathrm{VN}$ means trend and large $\mathrm{VN}$ means independence. 
Precisely, 
von Neumann~\cite{original_von_Neumann,earlier_von_Neumann_ratio_test} demonstrated that under the assumption of normality and independence, $E(\mathrm{VN})=2N/(N-1)$, so that $2$ can be taken as threshold value for large sample size. 
 In~\cite{animal}, an empirical distribution for $\mathrm{VN}$ was derived, appearing  
normal with mean $2$ for large $N$. 
Critical values of the left-tailed test of the Null Hypothesis of independence are derived in~\cite[Table 1]{animal}. 
%It is intuitively easy to see that ``small $\eta$" means ``existence of trend." 
%Precisely, 
%if the observations are drawn from a normal population, the theoretical distribution of $\eta$ as derived in~\cite{original_von_Neumann,earlier_von_Neumann_ratio_test}, 
%confronted with the experimental  value of $\mathrm{VN}$, could determine whether a trend exists. 

The problem is that, as observed in~\cite{Rank_von_Neumann}, this test is not robust against 
deviation from normality in the data. We therefore have to resort to a nonparametric test. 

\subsection{Rank von Neumann ratio test}

%In the von Neumann analysis, 
%independence is understood as independence of  time-successive observations $y_n$, $y_{n+1}$, ...  
%Along those lines, 
%we will mention the von Neumann ratio test~\cite{original_von_Neumann,earlier_von_Neumann_ratio_test}. 
%Here, however, there is no objective way to assign a ``time stamp" to the observations. 
%Here, however, the variables are rather categorical, 
%as there is no time variable that can be associated with the outcomes of the search algorithm.  

%The initial values are in the data set
%(results{X}.init_{bias,time}).

The nonparametric von Neumann test does not rely on numerical values of the observations, but on their ranking. 
%Even though the results cannot be ordered time-wise, it is nevertheless possible to 
The log sensitivity observations are ranked consistently with increasing error. 
To be specific, let $x_n$ be the error and $y_n$ the log-sensitivity with ``time-stamp" $n$ in bin $i$. 
Let $\pi_i : \{1,...,N_i\} \to \{1,...,N_i\}$ be the permutation of the set of $N_i$ labels such that 
$x_{\pi_i(n)}$ is nondecreasing.   
The $y_{\pi_i(n)}$ sequence may have an overall nondecreasing trend in case the classical limitation 
is violated, but it is not uniformly nondecreasing. We show randomness in the latter sequence. 
We define $r_n$ to be the rank of the observation $y_{\pi_i(n)}$ in a nondecreasing reordering of the data $\{y_n\}$. 
%the nonparametric {\it rank} version of the von Neumann test~\cite{Rank_von_Neumann} can be run.  
Inside the bin $i \in I$, the {\it rank} von Neumann ratio is 
\[ \mathrm{RVN}_i=\frac{\sum_{n=1}^{N_i-1}(r_n-r_{n+1})^2}{\sum_{n=1}^{N_i}(r_n-\bar{r}_i)^2}, 
\quad \left(\bar{r}_i=\frac{1}{N_i}\sum_{n=1}^{N_i} r_n \right).\]
The statistic of the $\mathrm{RVN}$ is approximately $\beta$-distributed, from which 
the left-tailed critical values of $\mathrm{RVN}_\alpha$ are cataloged in~\cite[Table 2]{Rank_von_Neumann}. 
%For example, for $N_i \approx 100$, $\mathrm{RVN}_{0.05}=1.67$. 
Naturally, in case the classical limitation is challenged, 
some correlation should be expected as the data $\{y_n\}$  is on the average increasing. 
For this reason, before running the von Neumann test, the ``trend" should be removed. 
%However, the rank test still reveals randomness. 

For example, consider a 7-ring with target spin 3, and 15 groups for a total of 1515 observations.
\begin{comment} 
the rank ratio test yields, 
%
\begin{eqnarray*}
\lefteqn{\mathrm{RVN}=}\\
&&1.9315  ~  2.4791  ~  1.8281 ~   2.2123  ~  1.6447 ~   2.2917  ~  1.5993  ~  1.6046  ~  1.7265 ~   1.2411\\
&&1.5564  ~  1.5317  ~  2.1542 ~   1.8349  ~  1.5842 
\end{eqnarray*}
%
\end{comment}
From~\cite[Table 2]{Rank_von_Neumann} of the thresholds of the $\beta$-function statistic, 
and observing that each group contains about 100 observations, 
any value $\geq 1.67$ would indicate randomness with 95\% confidence. 
%Clearly, randomness is present in most of the groups. 
The data $\{y_n\}_{n=1}^{N_i}$ is detrended 
with the {\tt detrend} function of Matlab, 
which removes the best overall linear fit 
of the trend from the data in bin $i$. Repeating this for all bins, 
the rank von Neumann ratio test for all 15 bins yields
\begin{eqnarray*}
\lefteqn{\mathrm{RVN}=}\\
&&2.1908 ~   1.8130 ~   2.1426  ~  1.7964  ~  1.7150  ~  2.0761  ~  1.3343  ~  1.9043  ~  1.7925  ~  1.2026\\
&&1.8994 ~   2.0175 ~   1.9932  ~  1.8010  ~  1.8657
\end{eqnarray*}
Observe that 13 out of 15 bins %(as opposed to 8 out of 15 in the nondetrended test) 
are beyond the threshold of 1.67, reinforcing our claim of randomness in the data.

\end{document}